\documentclass[11pt,reqno]{amsart}
\usepackage{amsfonts, amssymb, amscd}
\usepackage{amsmath}
\usepackage{verbatim}
\usepackage{amssymb}
\usepackage{mathrsfs}
\usepackage{graphicx}
\usepackage{cite}
\usepackage[all]{xy}
\usepackage{tikz}
\usepackage[toc,page]{appendix}
\usepackage{tikz-cd}

\usepackage{graphicx}
\usepackage{hyperref}
\usepackage{float}

\def\rk {{\operatorname{rk}}}
\def\C {{\mathbb C}}
\def\CC {{\mathbb C}}
\def\R {{\mathbb R}}
\def\Z {{\mathbb Z}}
\def\PP {{\mathbb P}}

\def\ii {{\rm i}}

\def\Vol {{\rm Vol}}
\allowdisplaybreaks[3]

\theoremstyle{definition}

\newtheorem{theorem}{Theorem}[section]
\newtheorem{lemma}[theorem]{Lemma}
\newtheorem{proposition}[theorem]{Proposition}
\newtheorem{definition}[theorem]{Definition}

\newtheorem{corollary}[theorem]{Corollary}
\newtheorem{remark}[theorem]{Remark}

\numberwithin{equation}{section}

\begin{document}
\title[Central charges in local mirror symmetry]{Central charges in local mirror symmetry via hypergeometric duality}

\begin{abstract}
We apply the better-behaved GKZ hypergeometric systems to study toric Calabi-Yau Deligne-Mumford stacks and their Hori-Vafa mirrors given by affine hypersurfaces in $(\C^*)^d$. We show that the integral structures of A-branes and B-branes coincide. This confirms a local version of a conjecture of Hosono and can be seen as a generalization of the Gamma conjecture for local mirror symmetry.


\end{abstract}

\author{Zengrui Han}
\address{Department of Mathematics\\
University of Maryland, College Park\\
College Park, MD, 20740} \email{zhan223@umd.edu}

\maketitle

\tableofcontents

\section{Introduction}\label{introduction}

\subsection{Toric Calabi-Yau Deligne-Mumford stacks and their Landau-Ginzburg mirrors}\label{notation}

Let $N$ be a lattice of rank $d+1$, and $\Delta$ be a lattice polytope of dimension $d$ which lies on a primitive hyperplane $\deg(-)=1$ in $N_{\R}:=N\otimes\R$ where $\deg: N\rightarrow \Z$ is a linear function. Let $C$ be the $(d+1)$-dimensional finite rational polyhedral cone in $N_{\R}$ over the polytope $\Delta$, i.e., $C:=\R_{\geq 0}(\Delta\oplus 1)$. This data encodes an affine toric variety with Gorenstein singularities $X=\operatorname{Spec}\CC[C^{\vee}\cap N^{\vee}]$, in the sense that its dualizing sheaf is trivial.

\smallskip

Let $\{v_i\}_{i=1}^n$ be a subset of lattice points in the polytope $\Delta$ that includes all the vertices. For each regular triangulation of the polytope $\Delta$ that involves some points in $\{v_i\}_{i=1}^n$, the corresponding subdivision $\Sigma$ of the cone $C$ defines a smooth toric Deligne-Mumford stack $\PP_{\mathbf{\Sigma}}$ in the sense of Borisov-Chen-Smith \cite{BCS}, where $\mathbf{\Sigma}:=(\Sigma,\{v_i\})$ is the stacky fan whose distinguished lattice points on the rays are chosen to be the primitive ray generators $v_i$'s, providing a crepant resolution of the affine Gorenstein toric variety $X$. The toric stack $\PP_{\mathbf{\Sigma}}$ is Calabi-Yau in the sense that its canonical class is trivial. 

\smallskip

Consider the hypersurface defined as $Z_f:=(f=0)\subseteq (\CC^*)^d$ where $f$ is a Laurent polynomial defined as
\begin{align*}
	f(z_1,\cdots,z_d)=\sum_{i=1}^n x_i z^{\bar{v_i}}=\sum_{i=1}^n x_i\cdot z_1^{v_{i1}}\cdots z_d^{v_{id}}
\end{align*}
Here, we denote by $\bar{v_i}$ the $d$-dimensional vector obtained from $v_i$ by deleting the last coordinate $1$. In this paper, we adopt the idea that the family of hypersurfaces $\{Z_f\}$ should be viewed as the mirror family of the smooth toric Calabi-Yau Deligne Mumford stack $\PP_{\mathbf{\Sigma}}$. Note that this is a generalization of the Hori-Vafa mirror \cite{Hori-Vafa}, where each cone in $\Sigma$ is assumed to be unimodular, thus the corresponding toric stack is a smooth variety.

\smallskip

According to Kontsevich's homological mirror symmetry conjecture, the bounded derived category $D^b(\PP_{\mathbf{\Sigma}})$ of coherent sheaves (category of B-branes) on the toric stack $\PP_{\mathbf{\Sigma}}$ should be equivalent to some appropriately defined Fukaya-Seidal category of the Landau-Ginzburg model\footnote{It is suggested by Iritani to the author that it might be more natural to consider the Landau-Ginzburg model ($((\C^*)^{d+1},z_{d+1}\cdot f)$) of one-dimensional higher and its associated Fukaya-Seidel category. The period integrals considered in this paper could be obtained from the oscillatory integrals of this new LG model. It would be interesting to investigate this direction.} $((\C^*)^d,f)$ (category of A-branes):
\begin{align}\label{HMS}
	D^b(\PP_{\mathbf{\Sigma}})\xrightarrow{\sim} FS((\C^*)^d,f)
\end{align}

\smallskip

While the mathematical meaning of the former is clear, the definition of the latter seems to be subtle. However, in any sense, the objects in a reasonably defined Fukaya-type category of $((\C^*)^d,f)$ should be Lagrangian submanifolds of the torus $(\CC^*)^d$ which satisfies certain admissibility condition that take the potential $f$ into account.

\smallskip

Another mathematical formulation of mirror symmetry is provided by the SYZ conjecture. Roughly speaking, it predicts that two spaces that are mirror to each other should arise from dual special Lagrangian torus fibrations over some base space. In our setting, the fibration on the mirror side is given by the logarithmic map
\begin{align}\label{fibration}
	\operatorname{Log}:(\CC^*)^d\rightarrow \R^d, \ (z_1,\cdots,z_d)\mapsto (\log|z_1|,\cdots,\log|z_d|)
\end{align}

\smallskip

Furthermore, the SYZ conjecture predicts that under the equivalence \eqref{HMS}, the mirror cycles of line bundles (resp. skyscraper sheaves) on $\PP_{\mathbf{\Sigma}}$ should be Lagrangian sections (fibers) of the fibration \eqref{fibration}. In particular, when the coefficients of the Laurent polynomial $f$ are all positive real numbers, the mirror cycle of the structure sheaf $\mathcal{O}_{\PP_{\mathbf{\Sigma}}}$ should be the positive real section $(\R_{>0})^d\subseteq (\CC^*)^d$, and those of arbitrary line bundles corresponding to torus-invariant divisors $\sum a_i D_i$ are obtained by parallel transport of $(\R_{>0})^d$.

\subsection{Central charges as solutions to better-behaved GKZ systems}\label{intro-central charges}

In this paper, instead of studying the HMS statement of the Hori-Vafa mirror pair in such generality, we focus on the more classical version of mirror symmetry, namely the equality between the \textit{central charges} of the A-branes and the B-branes.

\smallskip

The notion of the central charge appears in physics literatures and has played an important role in Douglas's $\Pi$-stability \cite{Douglas} and its mathematical formalism \cite{Bridgeland} due to Bridgeland. By definition the central charge of A-branes (B-branes, respectively) is a linear function defined on the Grothendieck group of the  derived category of A-branes (B-branes, respectively). 

\smallskip

In \cite{Hosono} Hosono defined the central charges for Calabi-Yau complete intersections in toric varieties (Batyrev-Borisov mirrors) in terms of period integrals and hypergeometric series respectively. The spaces of A- and B-branes are defined as certain homology groups associated to the pair $((\CC^*)^d, Z_f)$ and Grothendieck groups of the Calabi-Yau complete intersection $X$, which comes naturally with integral structures by definition. In \cite[Conjecture 2.2]{Hosono} Hosono conjectured that these two integral structures coincide under certain homological mirror symmetry correspondence. A local version of this conjecture can be made in a similar manner (see \cite[Conjecture 6.3]{Hosono}). In this paper, we focus on the local mirror symmetry setting, i.e., Hori-Vafa mirrors defined in the previous section.


\smallskip

It is also observed that both of the A-brane and B-brane central charges can be viewed as solutions to certain GKZ hypergeometric systems, therefore GKZ systems play an important role in the study of local mirror symmetry. However, we believe the \textit{better-behaved GKZ systems}, introduced by Borisov-Horja \cite{BH} and studied in a series of papers \cite{BHconj, BHW, BHanDuality, bbGKZ_analytic_continuation}, are more suitable for this consideration:
\begin{definition}\label{defGKZ}
	Consider the system of partial differential equations on the collection of functions $\{\Phi_c(x_1,\ldots,x_n)\}$ in complex variables $x_1,\ldots, x_n$, indexed by the lattice points in $C$:
	\begin{align*}
		\partial_i\Phi_c=\Phi_{c+v_i},\quad \sum_{i=1}^n\langle\mu,v_i\rangle x_i\partial_i\Phi_c+\langle\mu,c\rangle\Phi_c=0
	\end{align*}
for all $\mu\in N^{\vee}$, $c\in C$ and $i=1,\ldots,n$. We denote this system by $\mathrm{bbGKZ}(C,0)$. Similarly by considering lattice points in the interior $C^{\circ}$ only, we can define $\mathrm{bbGKZ}(C^{\circ},0)$.
\end{definition}

\begin{remark}\label{A-model interpretation}
	This definition is essentially a B-model interpretation. In terms of quantum cohomology (A-model), the better-behaved GKZ systems $\mathrm{bbGKZ}(C,0)$ and $\mathrm{bbGKZ}(C^{\circ},0)$ corresponds to the ordinary and compactly supported quantum cohomology $D$-module $QH^*(\PP_{\mathbf{\Sigma}})$ and $QH^*_c(\PP_{\mathbf{\Sigma}})$ of the toric stack $\PP_{\mathbf{\Sigma}}$ respectively. For details, see e.g. \cite{I-1}.
\end{remark}

\smallskip

It is proved in \cite{BHanDuality} that the two systems $\mathrm{bbGKZ}(C,0)$ and $\mathrm{bbGKZ}(C^{\circ},0)$ are dual to each other by constructing an explicit non-degenerate pairing between the solution spaces. The key observation in this paper is that the A-brane central charges (i.e., the period integrals) are naturally solutions to the \textit{dual} system $\mathrm{bbGKZ}(C^{\circ},0)$ rather than $\mathrm{bbGKZ}(C,0)$ (which is equivalent to the usual GKZ systems in normal cases). Motivated by this, we introduce our version of the A-brane central charges, which is a slight modification of the original definition of Hosono.

\smallskip

\begin{definition}
	For each lattice point $c$ in the interior $C^{\circ}$ of the cone $C$, we define the following holomorphic form
	\begin{align*}
		\omega_c:=(-1)^{\deg c-1}(\deg c-1)!\frac{z^c}{f^{\deg c}}\frac{dz_1}{z_1}\wedge\cdots\wedge\frac{dz_d}{z_d}
	\end{align*}
	on the complement $(\CC^*)^d\backslash Z_f$, where $f=\sum x_i z^{v_i}$. The \textit{A-brane central charge} associated to a Lagrangian submanifold $L$ (with certain admissibility conditions) is defined to be the collection of period integrals
	\begin{align*}
		Z^{A,L}(x)=(Z^{A,L}_c(x))_{c\in C^{\circ}}=\left(\frac{(-1)^d}{(2\pi{\rm i})^{d+1}}\int_{L}\omega_c\right)_{c\in C^{\circ}}.
	\end{align*}
	where each $\int_{L}\omega_c$ is viewed as a holomorphic function with the coefficients $x_i$ of $f$ as variables.
\end{definition}

\smallskip

Following the arguments of Batyrev \cite{Batyrev} and Borisov-Horja \cite{BH}, it is easy to verify that the A-brane central charges are solutions to the dual system $\mathrm{bbGKZ}(C^{\circ},0)$, see Proposition \ref{solution to bbGKZ}. The A-brane central charge map 
\begin{align*}
	Z^A:\ H_d((\CC^*)^d\backslash Z_f,\Z) \rightarrow \mathrm{Sol}(\mathrm{bbGKZ}(C^{\circ},0))
\end{align*}
defined here gives an integral structure on the solution spaces of the bbGKZ systems, by considering the image of the lattice $H_d((\CC^*)^d\backslash Z_f,\Z)$.

\smallskip

On the other hand, we use the hypergeometric series $\Gamma^{\circ}$ with values in the \textit{cohomology with compact support} as our definition of B-brane central charges. Note that in the original paper \cite{Hosono} Hosono used the usual hypergeometric series $\Gamma$ (that takes values in the usual cohomology space) as central charges on the B-side. For precise definitions of $\Gamma$ and $\Gamma^{\circ}$, see Definition \ref{definition of gamma series}.


\smallskip

\begin{definition}
	Let $N$, $C$, and $\mathbf{\Sigma}=(\Sigma,\{v_i\})$ be as previous. For any class $\mathcal{E}\in K_0(\PP_{\mathbf{\Sigma}})$ in the K-theory of $\PP_{\mathbf{\Sigma}}$, we define its \textit{B-brane central charge} by
	\begin{align*}
		Z^{B,\mathcal{E}}(x)=(Z^{B,\mathcal{E}}_c(x))_{c\in C^{\circ}}=(\chi(ch(\mathcal{E}),-)\circ \Gamma^{\circ}_c)_{c\in C^{\circ}}
	\end{align*}
\end{definition}

\smallskip

Similar to the A-side, the B-brane central charge map
\begin{align*}
	Z^B:\ K_0(\PP_{\mathbf{\Sigma}}) \rightarrow \mathrm{Sol}(\mathrm{bbGKZ}(C^{\circ},0))
\end{align*}
defines another integral structure on the solution space of bbGKZ systems by the natural $\Z$-structure on the Grothendieck group $K_0(\PP_{\mathbf{\Sigma}})$.

\smallskip

The main result of this paper is the matching of the integral structures on the A- and B- sides.


\smallskip

\begin{theorem}\label{main theorem}
	The A- and B-brane integral structures of the Hori-Vafa mirrors, defined by $H_{d}\left((\CC^*)^{d}\backslash Z_f,\Z\right)$ and 
$K_0(\PP_{\Sigma},\Z)$ respectively, coincide.
\end{theorem}




\subsection{Relations to Gamma conjecture and hypergeometric duality}

In \cite{AGIS}, Abouzaid, Ganatra, Iritani and Sheridan proposed a new approach to study the asymptotic behavior of period integrals in the setting of Batyrev mirror pairs (i.e., Calabi-Yau hypersurfaces in Fano toric varieties). Roughly speaking, the idea is to cut the cycles into pieces according to the tropical geometry of the Laurent polynomial $f$, then approximate each piece by the volume of certain polytopes which could again be related to the gamma classes on the mirror side. They utilized this method to give an alternative proof of the mirror-symmetric Gamma conjecture for the Batyrev mirror pairs.

\smallskip

The first step to prove our main result Theorem \ref{main theorem} is to apply their tropical method to our local mirror symmetry setting. More precisely, in \S\ref{asymptotic} we prove the following asymptotic formula for the period integral of the affine hypersurface $Z_f$ in $(\CC^*)^d$ corresponding to the positive real locus $\R_{>0}^d$.

\smallskip

\begin{theorem}[=Theorem \ref{asymptotics thm}]
	Fix a triangulation $\Sigma$ of the cone $C$, and denote the corresponding convex piecewise linear function by $\psi$. Then the asymptotic behavior of $\Psi_c(t^{-\psi(v_1)},...)$ when $t\rightarrow +\infty$ is given by
	\begin{align*}
		Z^{A,\R_{>0}^d}_c(t^{-\psi(v_1)},...)=&t^{\psi(c)}\frac{(-1)^{\rk N-\deg c}}{(2\pi{\rm i})^{|\sigma(c)|}|\operatorname{Box}(\sigma(\gamma))|}\cdot\int_{\gamma(c)} t^{\omega}\cdot \hat{\Gamma}_{\gamma(c)} F_{I_c}+o(t^{\psi(c)})
	\end{align*}
	For precise definitions of the symbols, see \S\ref{asymptotic}.
\end{theorem}

\smallskip

This could be viewed as the Gamma conjecture for local mirror symmetry. Our computation could be seen as a generalization of the one in \cite{AGIS} because they assumed the polytopes $\Delta$ are reflexive, while ours are arbitrary. Moreover, we compute the asymptotics of period integrals corresponding to arbitrary interior points in the cone $C$, while they essentially dealt with the one corresponding to the unique interior point of the (reflexive) polytope of degree equals to 1.

\smallskip

The second essential ingredient of the proof is the hypergeometric duality established by Borisov and the author in \cite{BHanDuality}. More precisely, an explicit formula for the pairing between the solution spaces to $\mathrm{bbGKZ}(C,0)$ and $\mathrm{bbGKZ}(C^{\circ},0)$ is constructed: 
\begin{align*}
		\langle \Phi,\Psi\rangle = \sum_{c,d,I}\xi_{c,d,I}\Vol_I\left(\prod_{i\in I}x_i\right)\Phi_c\Psi_d
\end{align*}
where $\Phi$ and $\Psi$ are solutions to $\mathrm{bbGKZ}(C,0)$ and $\mathrm{bbGKZ}(C^{\circ},0)$ respectively, and $\xi_{c,d,I}$'s are certain constant coefficients which only depends on the combinatorics of the cone. The non-degeneracy of this pairing allows us to lift the equality of leading terms to a genuine equality between central charges on the two sides. An immediate benefit of this is the following. In \cite{AGIS}, the authors has to construct explicitly tropical approximations of cycles mirror to arbitrary torus-invariant divisors $\sum a_iD_i$ and repeat the computation of asymptotics for them (because generally asymptotic expansion does not commute with taking leading terms). While in principle, we can do the same thing in our setting, it is no longer necessary since the equality on the level of functions themselves (rather than just the leading terms) allows us to conclude simply by looking at the monodromy on both sides, see Corollary \ref{monodromy}.

\subsection{Comparison with related results}

Let us discuss the relationship between the main result of this paper and other related works. In \cite{I-1} and \cite{I-2}, Iritani studied exponential period integrals (also known as oscillatory integrals) mirror to nef complete intersections in compact toric DM stacks using different methods. The results in these papers are more general than the ones in \cite{AGIS} in the sense that the period integrals are identified with explicit hypergeometric series, not just  the leading asymptotics. The period integrals considered in this paper could be presumably evaluated by using similar methods.

\smallskip

The computation in this paper is more general than that in \cite{I-2} in the sense that the polytope is not assumed to be reflexive. Therefore the hypersurface $Z_f$ we considered is not necessarily birational to a Calabi-Yau variety, and could have geometric genus larger than 1.

\smallskip

It is also pointed out by Hiroshi Iritani to the author that a recent paper of Yamamoto \cite{Yamamoto} is related to this work. It seems plausible that some of Yamamoto's results would follow from the result of this paper. Recent work \cite{FLYZ} of Fang, Liu, Yu and Zong studied remodeling conjecture with descendants, which is also closely related to Hosono's conjecture for toric Calabi-Yau 3-folds.

\smallskip

As is explained in Remark \ref{A-model interpretation}, the better-behaved systems could be seen as ordinary or compactly supported quantum cohomology $D$-modules. The natural central charges for these systems are given by compactly supported branes for the usual system and possibly non-compact branes for the dual system. The focus of this paper is on the latter. The former, which is closer to the original version of Hosono's conjecture, has been studied in \cite{Wang} using tropical geometry approach based on the foundational work of Ruddat-Siebert \cite{RS}, and in \cite{I-3} using methods similar to \cite{I-1} and \cite{I-2}. In \cite{I-3}, only local Calabi-Yau varieties arising as the total space of a canonical line bundle were investigated, which is a special case of the general toric stack $\PP_{\mathbf{\Sigma}}$ considered in this paper.

\subsection{Future directions}

The construction of A-brane central charges and the main result of this paper are not completely satisfactory. More specifically, the space $H_d((\CC^*)^d\backslash Z_f,\CC)$ we used as the space of A-branes is not of the correct dimension, see Section \ref{issues} for a more detailed discussion. Hopefully, a better understanding of the Grothendieck group of (certain version of) Fukaya-Seidel category of the pair $((\CC^*)^d, f)$ would yield the ``correct'' space of A-branes for the Hori-Vafa mirror pairs we considered in this paper.

\subsection{Organization of the paper}
The paper is organized as follows. In section \ref{preliminary} we recall basic definitions and properties of smooth toric Deligne-Mumford stacks and their orbifold cohomology and K-groups. In section \ref{central charges}, we formulate the precise definitions of A-brane and B-brane central charges in our setting, and interpret them as solutions to certain better-behaved GKZ systems. In section \ref{asymptotic}, we compute the asymptotic behavior of the period integrals (A-brane central charges) when approaching a large radius limit point corresponding to a regular triangulation $\Sigma$. In section \ref{equality}, we apply the hypergeometric duality to prove the equality of A-brane and B-brane central charges. Appendices \ref{volume} and \ref{beta function} contain auxiliary results that are used in the computation of section \ref{asymptotic}. 

\subsection{Acknowledgements}

The author would like to thank his advisor Lev Borisov for consistent support, helpful discussions and useful comments throughout the preparation of this paper. The author would like to thank Hiroshi Iritani and Chiu-Chu Melissa Liu for their interest in this work, warm encouragement, and valuable feedback on an earlier draft. The author thanks Mohammed Abouzaid, Bohan Fang, Helge Ruddat, Uli Walther and Shaozong Wang for useful conversations. The author thanks the anonymous referee, whose comments and suggestions greatly improved this paper, for the careful and thoughtful reading of the text.

\section{Toric stacks, orbifold cohomology and K-groups}\label{preliminary}

In this section, we review basic knowledge on smooth toric Deligne-Mumford stacks with an emphasize on their orbifold cohomology spaces and K-theory, and fix the notations that will be used throughout this paper. The main references are \cite{BCS, BHconj, BHanDuality}.

\subsection{Smooth toric DM stacks and the twisted sectors}

Following \cite{BCS}, a smooth toric Deligne-Mumford stack associated to a stacky fan $\mathbf{\Sigma}=(\Sigma,\{v_1,\cdots,v_n\})$ and its twisted sectors are defined as follows.

\smallskip

\begin{definition}
		Let $C$, $N$, and $\mathbf{\Sigma}=(\Sigma,\{v_1,\cdots,v_n\})$ be combinatorial datum defined in \S\ref{notation}. Consider the open subset $U$ of $\CC^n$ defined by
		\begin{align*}
			U=\left\{(z_1,\cdots,z_n)\in\CC^n:\ \{i:z_i=0\}\in\Sigma\right\}
		\end{align*}
		and a subgroup $G$ of $(\CC^*)^n$ defined by
		\begin{align*}
			G=\left\{(\lambda_1,\cdots,\lambda_n):\prod_{i=1}^n \lambda_i^{\langle m, v_i \rangle}=1,\forall m\in N^{\vee}\right\}
		\end{align*}
		The smooth toric Deligne-Mumford stack $\PP_{\mathbf{\Sigma}}$ associated to the stacky fan $\mathbf{\Sigma}$ is defined to be the stack quotient of $U$ by $G$. We denote its coarse moduli space (i.e., the toric variety associated to the ordinary fan $\Sigma$) by $\PP_{\Sigma}$.
\end{definition}

\smallskip

\begin{definition}
	For each cone $\sigma\in\Sigma$ we define $\operatorname{Box}(\sigma)$ to be the set of lattice points $\gamma$ which can be written as $\gamma=\sum_{i\in\sigma}\gamma_i v_i$ with $0\leq\gamma_i<1$. We denote the union of all $\operatorname{Box}(\sigma)$ by $\operatorname{Box}(\Sigma)$. To each element $\gamma\in\operatorname{Box}(\Sigma)$ we associate a twisted sector\footnote{Abusing the notation slightly, we will denote the closed substack $\PP_{\mathbf{\Sigma}/\sigma(\gamma)}$ by $\gamma$ when there is no confusion.} $\PP_{\mathbf{\Sigma}/\sigma(\gamma)}$ of $\PP_{\mathbf{\Sigma}}$ corresponding to the minimal cone $\sigma(\gamma)$ in $\Sigma$ containing $\gamma$, which is the closed substack of $\PP_{\mathbf{\Sigma}}$ defined by the quotient fan $\Sigma/\sigma(\gamma)$. We define the dual of a twisted sector $\gamma=\sum\gamma_i v_i$ by
\begin{align*}
	\gamma^{\vee}=\sum_{\gamma_i\not=0}(1-\gamma_i)v_i.
\end{align*}
or equivalently, the unique element in $\operatorname{Box}(\sigma(\gamma))$ that satisfies
\begin{align*}
	\gamma^{\vee}=-\gamma\hskip -3pt \mod \sum_{i\in\sigma}\Z v_i
\end{align*}
The inertia stack of $\PP_{\mathbf{\Sigma}}$ is the disjoint union of all the twisted sectors.
\end{definition}

\subsection{Orbifold cohomology and K-groups of toric stacks}

Due to the non-compactness of $\PP_{\mathbf{\Sigma}}$, there are two types of orbifold cohomology theory associated to it, namely the usual orbifold cohomology and the orbifold cohomology with compact support. The following results are proved in \cite{BHconj}.

\smallskip

\begin{proposition}\label{coh}
As usual, $\operatorname{Star}(\sigma(\gamma))$ denotes the set of cones in $\Sigma$ that contain $\sigma(\gamma)$.
	Cohomology space $H_{\gamma}$ of the twisted sector $\gamma$ is naturally isomorphic to the quotient of the polynomial ring $\CC[D_i:i\in\operatorname{Star}(\sigma(\gamma))\backslash\sigma(\gamma)]$ by the ideal generated by the relations
	\begin{align*}
		\prod_{j\in J}D_j,\ \text{for }J\not\in\operatorname{Star}(\sigma(\gamma))
	\end{align*}
	and
	\begin{align*}
		\sum_{i\in\operatorname{Star}(\sigma(\gamma))\backslash\sigma(\gamma)}\mu(v_i)D_i,\ \text{for } \mu\in\operatorname{Ann}(v_i,i\in\sigma(\gamma)).
	\end{align*}
	There is a $\CC[D_1,\ldots, D_n]$-module structure on $H_{\gamma}$ defined by declaring $D_i=0$ for $i\not\in \operatorname{Star}(\sigma(\gamma))$ and solving (uniquely) for 
$D_i, i\in \sigma(\gamma)$ to satisfy the linear relations $\sum_{i=1}^n \mu(v_i)D_i=0$ for all $\mu\in N^\vee$.
\end{proposition}

\smallskip

\begin{proposition}\label{cohcmp}
	Cohomology space with compact support $H_{\gamma}^c$ is generated by $F_I$ for $I\in\operatorname{Star}(\sigma(\gamma))$ such that $\sigma_I^{\circ}\subseteq C^{\circ}$ with relations
	\begin{align*}
		&\quad D_i F_I-F_{I\cup\{i\}}\text{ for }i\not\in I,I\cup\{i\}\in\operatorname{Star}(\sigma(\gamma))\\
		&\text{and }D_i F_I\text{ for }i\not\in I,I\cup\{i\}\not\in\operatorname{Star}(\sigma(\gamma))
	\end{align*}
as a module over $H_{\gamma}$.
\end{proposition}

\smallskip

There is a natural integration map $\int$ defined on each cohomology space $H_{\gamma}^c$ with compact support.

\smallskip

\begin{proposition}
	There exists a unique linear function $\int_{\gamma}:H_{\gamma}^c\rightarrow\C$ that takes values $\frac{1}{\Vol_{\overline{I}}}$ on each generator $F_I$ with $|I|=d+1-|\sigma(\gamma)|$ (i.e., of highest degree), where $\Vol_{\overline{I}}$ denotes the volume of the cone $\overline{\sigma_{I}}$ in the quotient fan $\Sigma/\sigma(\gamma)$. Moreover, it takes value zero on all elements 
of $H_{\gamma}^c$ of lower degree.
\end{proposition}
\begin{proof}
	See \cite[Proposition 2.6]{BHconj}.
\end{proof}

\smallskip

\begin{definition}
	The orbifold cohomology $H_{\mathrm{orb}}^*(\PP_{\mathbf{\Sigma}})$ of the smooth toric DM stack $\PP_{\mathbf{\Sigma}}$ is defined as the direct sum $\bigoplus_{\gamma} H_{\gamma}$ over all twisted sectors. Similarly, the orbifold cohomology with compact support $H_{\mathrm{orb},c}^*(\PP_{\mathbf{\Sigma}})$ is defined as $\bigoplus_{\gamma} H_{\gamma}^c$. We denote by $1_\gamma$ the generator of $H_\gamma$.
\end{definition}

\smallskip

\begin{remark}
	There is an involution map $*$ on the orbifold cohomology $H_{\mathrm{orb}}^*(\PP_{\mathbf{\Sigma}})$ that maps $H_{\gamma}$ to $H_{\gamma^{\vee}}$, defined by $(1_{\gamma})^*=1_{\gamma^{\vee}}$ and $(D_i)^*=-D_i$.
\end{remark}

\smallskip

There is a special type of characteristic classes of the smooth toric DM stack $\PP_{\mathbf{\Sigma}}$, called \textit{Gamma classes}, which play an essential role in the computation of this paper. A similar definition has been introduced in \cite{I-1} for general smooth DM stacks. The version we used in this paper comes from \cite[Corollary 3.14]{BHanDuality}.

\smallskip

\begin{definition}\label{def-gamma class}
For each twisted sector $\gamma$ of $\PP_{\mathbf{\Sigma}}$, we define its Gamma class by
\begin{align*}
	\widehat{\Gamma}_{\gamma} = \prod_{i\in\sigma(\gamma)}\Gamma(\gamma_i+\frac{D_i}{2\pi {\rm i}})\prod_{i\in\operatorname{Star}(\sigma(\gamma))\backslash\sigma(\gamma)}\Gamma(1+\frac{D_i}{2\pi {\rm i}})
\end{align*}
which is a cohomology class in $H_{\gamma}^*$. We define the Gamma class of $\PP_{\mathbf{\Sigma}}$ to be the direct sum of Gamma classes of all of its twisted sectors.
\end{definition}

\smallskip

Next we look at the K-groups of the toric stack $\PP_{\mathbf{\Sigma}}$. Again, there are two types of K-groups, the ordinary K-group $K_0(\PP_{\mathbf{\Sigma}})$ and the compactly supported K-group $K_0^c(\PP_{\mathbf{\Sigma}})$. The former is the usual Grothendieck group of the bounded derived category $D^b(\PP_{\mathbf{\Sigma}})$. The latter is defined as the Grothendieck group of the triangulated category $D^c(\PP_{\mathbf{\Sigma}})$ of bounded complexes of coherent sheaves supported on the union of all compact toric divisors of $\PP_{\mathbf{\Sigma}}$. We have the following combinatorial descriptions.

\smallskip

\begin{proposition}
	Let $C$, $v_i$ and $\Sigma$ be as before. We denote the class of the line bundle $\mathcal{O}_{\PP_{\mathbf{\Sigma}}}(D_i)$ corresponding to the ray $v_i$ by $R_i$. Then the Grothendieck group $K_0(\PP_{\mathbf{\Sigma}})$ is isomorphic to the quotient of the ring $\Z[R_i^{\pm 1}]$ by the relations
	\begin{align*}
		\prod_{i=1}^n R_i^{\mu(v_i)}-1,\quad \forall\mu\in N^{\vee},\text{ and } \prod_{i\in I}(1-R_i),\quad \forall I\not\in\Sigma
	\end{align*}
	Furthermore, if we denote the class of the structure sheaf of the closed substack corresponding to a cone $\sigma_I$ by $G_I$, then $K_0^c(\PP_{\mathbf{\Sigma}})$ is a module over $K_0(\PP_{\mathbf{\Sigma}})$ generated by $G_I$ for all $I\in\Sigma$ with $\sigma_{I}$ being an interior cone, with the relations given by
	\begin{align*}
	(1-R_i^{-1})G_I=G_{I\cup\{i\}}\text{ if }I\cup\{i\}\in\Sigma \text{ and }0\text{ otherwise},\ \forall i.
    \end{align*}
\end{proposition}

\begin{proof}
	See \cite[Proposition 3.3, Definition 3.9]{BHconj}.
\end{proof}

\smallskip

The complexified K-groups $K_0(\PP_{\mathbf{\Sigma}})_{\CC}:=K_0(\PP_{\mathbf{\Sigma}})\otimes_{\Z}\CC$ and $K_0^c(\PP_{\mathbf{\Sigma}}):=K_0^c(\PP_{\mathbf{\Sigma}})\otimes_{\Z}\CC$ are related to the orbifold cohomology $H_{\mathrm{orb}}^*(\PP_{\mathbf{\Sigma}})$ and $H_{\mathrm{orb},c}^*(\PP_{\mathbf{\Sigma}})$ by the combinatorial Chern characters as follows.

\smallskip

\begin{proposition}
	There is a natural isomorphism
	\begin{align*}
		ch: K_0(\PP_{\mathbf{\Sigma}})_{\CC}\xrightarrow{\sim} H_{\mathrm{orb}}^*(\PP_{\mathbf{\Sigma}})=\bigoplus_{\gamma\in\operatorname{Box}(\Sigma)}H_\gamma
	\end{align*}
	defined by
	\begin{align*}
		ch_{\gamma}(R_i)=\begin{cases}
1,\quad & i\not\in\operatorname{Star}(\sigma(\gamma)) \\
e^{D_i},\quad &i\in\operatorname{Star}(\sigma(\gamma))\backslash\sigma(\gamma) \\
e^{2\pi{\rm i}\gamma_i}\prod_{j\not\in\sigma(\gamma)}ch_{\gamma}(R_j)^{\mu_i(v_j)},\quad & i\in\sigma(\gamma)
\end{cases}
	\end{align*}
	Similarly, there is a natural isomorphism
	\begin{align*}
		ch^c: K_0^c(\PP_{\mathbf{\Sigma}})_{\CC}\xrightarrow{\sim} H_{\mathrm{orb},c}^*(\PP_{\mathbf{\Sigma}})=\bigoplus_{\gamma\in\operatorname{Box}(\Sigma)}H_\gamma^c
	\end{align*}
	defined by
	\begin{align*}
		ch_{\gamma}^c(\prod_{i=1}^n R_i^{l_i}G_I)=\begin{cases}
0,\quad & I\not\subseteq \operatorname{Star}(\sigma(\gamma)) \\
\prod_{i=1}^n ch_{\gamma}(R_i^{l_i})F_I,\quad &I \subseteq \operatorname{Star}(\sigma(\gamma)) 
\end{cases}
	\end{align*}
\end{proposition}
\begin{proof}
	See \cite[Proposition 3.7, 3.11]{BHconj}.
\end{proof}

\subsection{Euler characteristic pairing}

There is a natural non-degenerate pairing $\chi(-,-)$ between $K_0(\PP_{\mathbf{\Sigma}})$ and $K_0^c(\PP_{\mathbf{\Sigma}})$ called {\it Euler characteristic pairing} defined as the alternative sum of the dimension of Ext groups. More precisely, let $\mathcal{F}^{\bullet}$ and $\mathcal{G}^{\bullet}$ be complexes in the derived categories $D^b(\PP_{\mathbf{\Sigma}})$ and $D^c(\PP_{\mathbf{\Sigma}})$ respectively, we define
\begin{align*}
	\chi(\mathcal{F}^{\bullet},\mathcal{G}^{\bullet})=\sum_{i=0}^{\infty}\dim \operatorname{Hom}_{D^b(\PP_{\mathbf{\Sigma}})}(\mathcal{F}^{\bullet},\mathcal{G}^{\bullet}[i]).
\end{align*}
In particular, if we take $\mathcal{F}^{\bullet}$ to be the structure sheaf $\mathcal{O}_{\PP_{\mathbf{\Sigma}}}$ and $\mathcal{G}^{\bullet}$ to be a coherent sheaf, then this definition recovers the usual Euler characteristic of coherent sheaves. Note that this pairing is defined over $\Z$.

\smallskip

On the other hand, the pairing (after extending to $\CC$-coefficients) could be translated to the orbifold cohomology spaces $H_{\mathrm{orb}}^*(\PP_{\mathbf{\Sigma}})$ and $H_{\mathrm{orb},c}^*(\PP_{\mathbf{\Sigma}})$ via the Chern character defined in the last subsection. To avoid cumbersome notation, we denote the pairing on the orbifold cohomology spaces and on the K-groups by the same symbol $\chi$.

\smallskip

\begin{proposition}\label{def-eulerpairing}
The Euler characteristic pairing
\begin{align*}
	\chi:H_{\mathrm{orb}}^*(\PP_{\mathbf{\Sigma}})\otimes H_{\mathrm{orb},c}^*(\PP_{\mathbf{\Sigma}})\rightarrow\C
\end{align*}
on the toric DM stack $\mathbb{P}_{\Sigma}$ is given by
	\begin{align*}
		\chi(a,b)=\chi(\oplus_{\gamma}a_{\gamma},\oplus_{\gamma}b_{\gamma})=\sum_{\gamma}\frac{1}{|\operatorname{Box}(\sigma(\gamma))|}\int_{\gamma^{\vee}}\operatorname{Td}(\gamma^{\vee})a_{\gamma}^*b_{\gamma^{\vee}}
	\end{align*}
	Here $\operatorname{Td}(\gamma)$ is the Todd class of the twisted sector $\gamma$, defined as
	\begin{align*}
		\operatorname{Td}(\gamma)=\frac{\prod_{i\in\operatorname{Star}\sigma(\gamma)\backslash\sigma(\gamma)}D_i}{\prod_{i\in\operatorname{Star}\sigma(\gamma)}(1-e^{-D_i})}.
	\end{align*}
\end{proposition}
\begin{proof}
	See \cite[Lemma 4.20]{BHW}.
\end{proof}

\smallskip

The following easy consequence will be used in Appendix \ref{volume}.

\smallskip

\begin{corollary}\label{formula of euler characteristic}
	The Euler characteristic of the sheaves represented by the class $v\in K_0^c(\PP_{\mathbf{\Sigma}})$ is given by
	\begin{align*}
		\chi(v) = \sum_{\gamma\in\operatorname{Box}(\Sigma)}\frac{1}{|\operatorname{Box}(\sigma(\gamma))|}\int_{\gamma}ch_{\gamma}^c(v)\operatorname{Td}(\gamma).
	\end{align*}
\end{corollary}

\section{Central charges as solutions to better-behaved GKZ systems}\label{central charges}

In this section, we give precise definitions of the A-brane and B-brane central charges for Hori-Vafa mirror pairs. Our definitions differ slightly from the ones in \cite{Hosono}. Along the lines we briefly recall the basic preliminaries needed to formulate the main results of the paper.

\subsection{A-brane central charges}

Denote the coordinates on the torus $(\C^*)^d$ by $z=(z_1,\cdots,z_d)$. Consider the Laurent polynomial $f=\sum_{i=1}^n x_i z^{\bar{v_i}}$. As in \S\ref{notation}, $\bar{v_i}$ denotes the $d$-dimensional vector obtained from $v_i$ by deleting the last coordinate $1$. On the other hand, for a lattice point $c\in N_{\R}$, we write $c=(\bar{c},\deg c)$ where $\deg c$ is the last coordinate of $c$ and $\bar{c}$ consists of the first $d$ coordinates.

\smallskip

\begin{definition}\label{A-brane central charge}
	For each lattice point $c$ in the interior $C^{\circ}$ of the cone $C$, we define the following holomorphic form
	\begin{align*}
		\omega_c:=(-1)^{\deg c-1}(\deg c-1)!\frac{z^{\bar{c}}}{f^{\deg c}}\frac{dz_1}{z_1}\wedge\cdots\wedge\frac{dz_d}{z_d}
	\end{align*}
	on the complement $(\CC^*)^d\backslash Z_f$, where $f=\sum x_i z^{\bar{v_i}}$. The \textit{A-brane central charge} associated to a Lagrangian submanifold $L$ (with certain admissibility conditions) is defined to be the collection of period integrals\footnote{ Note that the constant factor $\frac{(-1)^d}{(2\pi{\rm i})^{d+1}}$ plays no essential role.}
	\begin{align*}
		Z^{A,L}(x)=(Z^{A,L}_c(x))_{c\in C^{\circ}}=\left(\frac{(-1)^d}{(2\pi{\rm i})^{d+1}}\int_{L}\omega_c\right)_{c\in C^{\circ}}
	\end{align*}
	where each $\int_{L}\omega_c$ is viewed as a holomorphic function with the coefficients $x_i$ of $f$ as variables.
\end{definition}

\smallskip

Since we are mostly interested in the mirror cycles of line bundles in this paper, from now on we will assume $L$ to be the Lagrangian sections of the fibration $\pi:(\CC^*)^d\rightarrow \R^d$ defined by $z\mapsto\log|z|$, see the discussion prior to \S\ref{intro-central charges}. The following result explains the reason why it is more natural to think of the period integrals over Lagrangian sections as solutions to the dual system $\operatorname{bbGKZ}(C^{\circ},0)$ rather than the usual system $\operatorname{bbGKZ}(C,0)$. Similar results can be found in \cite{EulerMellinIntegral1} and \cite{EulerMellinIntegral2}. 

\smallskip

\begin{lemma}
	The period integral
	\begin{align*}
		\int_{\R_{\geq 0}^d}\frac{z^{\bar{c}}}{f(z)^{\deg c}}\frac{dz_1}{z_1}\wedge\cdots\wedge\frac{dz_d}{z_d}
	\end{align*}
	is absolutely convergent if and only if $c\in C^{\circ}$.
\end{lemma}
\begin{proof}
	We make the coordinate change $z_i = e^{y_i}$, then the integral becomes
	\begin{align*}
		\int_{\R^d}\frac{e^{\bar{c}\cdot y}}{f(e^{y})^{\deg c}}dy_1\wedge\cdots\wedge dy_d =\int_{\R^d}\frac{e^{\bar{c}\cdot y}}{(\sum_{j\in\Delta} x_j e^{\bar{v_j}\cdot y})^{\deg c}}dy_1\wedge\cdots\wedge dy_d 
	\end{align*}
	Now we divide the space $\R^d$ into cone regions according to the normal fan $\Sigma$ of the polytope $\Delta$. More precisely, we divide $\R^d$ as the union of the following cone regions
	\begin{align*}
		\{\sigma_{\bar{v_k}}:=-\left(\R_{\geq 0}(\Delta-\bar{v_k})\right)^{\vee}:\bar{v_k}\in\Delta\}
	\end{align*}
	note that this differs from the usual definition of the normal fan by a negative sign. Fix a cone region $\sigma_{\bar{v_k}}$, it's easy to see that over $\sigma_{\bar{v_k}}$ the dominant term in the denominator $\sum_{j\in\Delta} x_j e^{\bar{v_j}\cdot y}$ is exactly the monomial $x_k e^{\bar{v_k}\cdot y}$. Therefore it suffices to consider the absolute convergence of 
	\begin{align*}
		\int_{\sigma_{v_k}}e^{(\bar{c}-(\deg c)\bar{v_k})\cdot y}dy_1\wedge\cdots\wedge dy_d 
	\end{align*}
	which is again equivalent to the condition that
	\begin{align}\label{condition_convergence}
		(\bar{c}-(\deg c)\bar{v_k})\cdot y<0,\quad \forall \text{ ray generators }y \text{ of }\sigma_{\bar{v_k}} \text{ and } \forall \bar{v_k}\in\Delta
	\end{align}
  Recall that the polytope $\Delta$	could be defined as the intersection of half-spaces (i.e., the facet representation of $\Delta$):
	\begin{align*}
		\Delta = \bigcap_{F \text{ facet of } \Delta} (h_F\geq 0)
	\end{align*}
	where $h_F$ is the defining equation of the support hyperplane spanned by the facet $F$. It is then straightforward to observe that the condition \eqref{condition_convergence} is equivalent to that $\bar{c}/(\deg c)$ is an interior point of the polytope $\Delta$, which is again equivalent to that $c$ is in the interior of the cone $C$.
\end{proof}

\smallskip

Now we consider a general Lagrangian section $L$. Additional restrictions are required to ensure the absolute convergence of the period integral. We write the section $L:\R^d\rightarrow(\CC^*)^d$ as
\begin{align*}
	(y_1,\cdots,y_d)\mapsto (e^{y_1}\cdot e^{{\rm i}\theta_1(y_1,\cdots,y_d)},\cdots,e^{y_d}\cdot e^{{\rm i}\theta_d(y_1,\cdots,y_d)})
\end{align*}

\smallskip

\begin{proposition}\label{convergence_of_the_integral}
	Suppose $L:\R^d\rightarrow(\CC^*)^d$ is a section of the fibration $T\rightarrow \R^d$ such that $\det( I_d+\left(\frac{\partial\theta_i}{\partial y_j}\right)_{i,j} )$ is bounded, then the integral $\int_{L}\omega_c$ is absolutely convergent for all $c\in C^{\circ}$.
\end{proposition}
\begin{proof}
	Follows directly from the last lemma and the observation that $\det( I_d+\left(\frac{\partial\theta_i}{\partial y_j}\right)_{i,j} )$ is the determinant of the Jacobian of the change of variables.
\end{proof}

\smallskip

\begin{remark}
	The condition on the Lagrangian section $L$ we made here is a bit artificial. We are not aware of a more natural condition that ensures the absolute convergence of the period integrals. 
\end{remark}

\smallskip

\begin{proposition}\label{solution to bbGKZ}
	Suppose $\gamma$ satisfies the condition in Proposition \ref{convergence_of_the_integral}, then $\Psi=(\Psi_c)_{c\in C^{\circ}}$ where
	\begin{align*}
		\Psi_c(x_1,\cdots,x_n):=\int_{L}\omega_c
	\end{align*}
	gives a solution to the system $\operatorname{bbGKZ}(C^{\circ},0)$.
\end{proposition}
\begin{proof}
	The idea of the proof comes from Batyrev \cite{Batyrev} and Borisov-Horja \cite{BH}. However, since the cycles we are integrating over are non-compact, some additional cares must be taken.
	
	To prove the equation $\partial_i\Psi_c=\Psi_{c+v_i}$ for any $i$, note that we have
	\begin{align*}
		\partial_i\left( \frac{z^{\bar{c}}}{f(z)^{\deg c}} \right) = (-\deg c)\frac{z^{\bar{c}+\bar{v_i}}}{f(z)^{\deg (c+v_i)}}
	\end{align*}
	which gives $\partial_i\omega_c = \omega_{c+v_i}$, and the absolute convergence of the integral ensures that differentiation commutes with integration.
	
	To prove the second equation
	\begin{align*}
		 \sum_{i=1}^n\langle\mu,v_i\rangle x_i\partial_i\Psi_c+\langle\mu,c\rangle\Psi_c=0,\quad \forall \mu\in N^{\vee}
	\end{align*}
	we look at the standard basis $\mu_1,\cdots,\mu_{d+1}$ of $N^{\vee}$. For $1\leq k\leq d$, an elementary computation shows that $\sum_{i=1}^n\langle\mu_k,v_i\rangle x_i\omega_{c+v_i}+\langle\mu_k,c\rangle\omega_c$ is equal to
	\begin{align*}
		(-1)^{\deg c-1}(\deg c-1)!z_k\partial_{z_k}\left(\frac{z^{\bar{c}}}{f^{\deg c}}\right)\frac{dz_1}{z_1}\wedge\cdots\wedge\frac{dz_d}{z_d}
	\end{align*}
	Note that
	\begin{align*}
		z_k\partial_{z_k}\left(\frac{z^{\bar{c}}}{f^{\deg c}}\right)\frac{dz_1}{z_1}\wedge\cdots\wedge\frac{dz_d}{z_d}=d\left(\frac{z^{\bar{c}}}{f^{\deg c}}\frac{dz_1}{z_1}\wedge\cdots\wedge\widehat{\frac{dz_k}{z_k}} \wedge\cdots\frac{dz_d}{z_d}\right)
	\end{align*}
	Take a chain of compact subsets $B_1\subseteq B_2\subseteq\cdots\subseteq T$ such that $\cup B_m =T$ (e.g., take $B_m$ to be the box defined by $e^{-m}\leq|z_j|\leq e^m$). By Stokes' theorem the integration of $z_k\partial_{z_k}\left(\frac{z^{\bar{c}}}{f^{\deg c}}\right)$ over $B_i\cap L$ is equal to
	\begin{align*}
		\int_{\partial(B_m\cap L)}\frac{z^{\bar{c}}}{f^{\deg c}}\frac{dz_1}{z_1}\wedge\cdots\wedge\widehat{\frac{dz_k}{z_k}} \wedge\cdots\frac{dz_d}{z_d}
	\end{align*}
	which tends to 0 when $m\rightarrow +\infty$ due to the absolute convergence of 
	\begin{align*}
		\int_{L}\frac{z^{\bar{c}}}{f^{\deg c}}\frac{dz_1}{z_1}\wedge\cdots\wedge\frac{dz_d}{z_d}.
	\end{align*}
	On the other hand, by dominated convergence theorem the sequence of integrals converges to the integration of $z_k\partial_{z_k}\left(\frac{z^{\bar{c}}}{f^{\deg c}}\right)$ over $L$. This finishes the proof of the case when $1\leq k\leq d$.
	Finally, if $k=d+1$, i.e., $\mu_k=\deg$, then an elementary computation shows that $\sum_{i=1}^n\langle\mu_k,v_i\rangle x_i\omega_{c+v_i}+\langle\mu_k,c\rangle\omega_c$ is zero.
\end{proof}

\smallskip

It is clear that a Lagrangian section $L$ satisfying the condition in Proposition \ref{convergence_of_the_integral} defines an integral homology class in the group $H_d((\CC^*)^d\backslash Z_f,\Z)$, therefore the map
\begin{align*}
	Z^A:\ H_d((\CC^*)^d\backslash Z_f,\Z) \rightarrow \mathrm{Sol}(\mathrm{bbGKZ}(C^{\circ},0))
\end{align*}
defines an integral structure on the space of solutions to the bbGKZ system associated to $C^{\circ}$. Unlike the B-brane integral structure that will be defined in the next subsection, which is only locally defined near the large radius limits, this is a globally defined integral structure over the whole stringy K\"ahler moduli space.

\subsection{B-brane central charges}

In this subsection we define central charges on the B-brane, i.e., the toric Deligne-Mumford stack $\PP_{\mathbf{\Sigma}}$, in terms of certain cohomology-valued Gamma series. We will be using the same notation as \S\ref{notation}. 

\smallskip

\begin{definition}\label{definition of gamma series}
We define the cohomology-valued Gamma series $\Gamma$ and $\Gamma^{\circ}$ as
\begin{align}
\Gamma_c(x_1,\ldots,x_n) =\bigoplus_{\gamma} \sum_{l\in L_{c,\gamma}} \prod_{i=1}^n 
\frac {x_i^{l_i+\frac {D_i}{2\pi \ii }}}{\Gamma(1+l_i+\frac {D_i}{2\pi \ii })}
\end{align}
for lattice point $c\in C$ and
\begin{align*}
\Gamma_c^{\circ}(x_1,\ldots,x_n) =\bigoplus_{\gamma} \sum_{l\in L_{c,\gamma}} \prod_{i=1}^n 
\frac {x_i^{l_i+\frac {D_i}{2\pi \ii }}}{\Gamma(1+l_i+\frac {D_i}{2\pi \ii })}\left(\prod_{i\in\sigma}D_i^{-1}\right)F_{\sigma}
\end{align*}
for lattice point $c\in C^{\circ}$, where both direct sums are taken over twisted sectors $\gamma = \sum_{j\in \sigma(\gamma)} \gamma_j v_j$ and the set $L_{c,\gamma}$ is the set of solutions to $\sum_{i=1}^n l_i v_i = -c$ with $l_i-\gamma_i \in \Z$ for all $i$, and $\sigma$ is the set of $i$ with $l_i\in\Z_{<0}$. The numerator is defined by picking a branch of $\log(x_i)$.
\end{definition}

\smallskip

It is proved in \cite{BHanDuality} that these series converge absolutely and uniformly on compacts in a neighborhood of the large radius limit point corresponding to the triangulation $\Sigma$. After composing them with linear functions on the orbifold cohomology spaces, we get holomorphic functions with values in $\CC$. It is proved that all solutions to the systems $\operatorname{bbGKZ}(C,0)$ and $\operatorname{bbGKZ}(C^{\circ},0)$ are obtained by composing $\Gamma$ and $\Gamma^{\circ}$ with linear functions on $H_{\mathrm{orb}}^*(\PP_{\mathbf{\Sigma}})$ and $H_{\mathrm{orb},c}^*(\PP_{\mathbf{\Sigma}})$.

\smallskip

Now we are able to define our version of B-brane central charge.

\smallskip

\begin{definition}\label{B-brane central charge}
	Let $N$, $C$, and $\mathbf{\Sigma}=(\Sigma,\{v_i\})$ be as previous. For any class $\mathcal{E}\in K_0(\PP_{\mathbf{\Sigma}})$ in the K-theory of $\PP_{\mathbf{\Sigma}}$, we define its \textit{B-brane central charge} by
	\begin{align*}
		Z^{B,\mathcal{E}}(x)=(Z^{B,\mathcal{E}}_c(x))_{c\in C^{\circ}}=(\chi(ch(\mathcal{E}),-)\circ \Gamma^{\circ}_c)_{c\in C^{\circ}}.
	\end{align*}
\end{definition}

\smallskip

Similar to the A-brane side, the map
\begin{align*}
	Z^B:\ K_0(\PP_{\mathbf{\Sigma}}) \rightarrow \mathrm{Sol}(\mathrm{bbGKZ}(C^{\circ},0))
\end{align*}
defines a B-brane integral structure on the space of solutions to the bbGKZ system associated to $C^{\circ}$. Note that this integral structure is only defined locally near the large radius limit corresponding to the triangulation $\Sigma$.

\subsection{Comments}\label{issues}

Here we make some comments on the integral structures of the bbGKZ systems defined by A- and B-brane central charges in this section.

According to the results of a series of papers of Borisov and Horja \cite{BHconj,BH,MellinBarnes}, the dimension of the complexified K-group $K_0(\PP_{\mathbf{\Sigma}})_{\CC}$ is equal to the normalized volume of the polytope $\Delta$, which is again equal to the dimension of the space of solutions to the bbGKZ system $\operatorname{bbGKZ}(C^{\circ},0)$. In fact, the B-brane central charge map $Z^B$ defined in the previous section is an isomorphism of vector spaces.

However, the A-brane central charge map
\begin{align*}
	Z^A: H_d((\CC^*)^d\backslash Z_f, \CC)\rightarrow \mathrm{Sol}(\mathrm{bbGKZ}(C^{\circ},0))
\end{align*}
is not an isomorphism for a simple reason -- the dimension of the homology space is strictly larger\footnote{According to Batyrev \cite{Batyrev}, the dimension of the $\CC$-coefficient homology space is equal to $\operatorname{vol}(\Delta) + d$.} than the dimension of the solution space. Therefore the construction is not completely satisfactory. We nevertheless admit this drawback and take the image of the integral homology group $H_d((\CC^*)^d\backslash Z_f, \Z)$ under the central charge map $Z^A$ as our definition of the integral structure of the bbGKZ systems on the A-brane side. The main result of this paper nevertheless says that we can find a sublattice generated by explicitly defined Lagrangian sections of $H_d((\CC^*)^d\backslash Z_f, \Z)$ of the correct rank that matches the integral structure defined by $K_0(\PP_{\mathbf{\Sigma}})$ on the B-side.

\section{Asymptotic behavior of period integrals via tropical geometry}\label{asymptotic}

The goal of this section is to analyze the asymptotic behavior of the A-brane central charge associated to the real positive locus $(\R_{> 0})^d$.

\smallskip

 To begin, we set up the notations that will be used throughout this section. We will be using the same notations as in \S\ref{central charges}. Additionally, let $\Sigma$ be a regular triangulation of the cone $C$, and denote the corresponding convex piecewise linear function by $\psi$. For each lattice point $c$ in the interior $C^{\circ}$, we denote the minimal cone in $\Sigma$ containing $c$ by $\sigma(c)$, and write $c=\sum_{i\in\sigma(c)}c_iv_i$. Then there is a unique twisted sector $\gamma(c)\in\operatorname{Box}(\Sigma)$ given by $\sum_{i\in\sigma(c)}\{c_i\}v_i$, where $\{c_i\}$ denotes the fractional part of $c_i$. Finally, we denote the set of indices $i$ such that $c_i=1$ by $I_c$. Note that $\sigma(c)$ is a disjoint union of $I_c$ and $\sigma(\gamma(c))$.
 
 \smallskip

The main result of this section is an asymptotic formula for the A-brane central charges of the positive real Lagrangian $\R_{>0}^d$ when the parameter approaches the large radius limit corresponding to the triangulation $\Sigma$, which should correspond to the structure sheaf $\mathcal{O}_{\PP_{\mathbf{\Sigma}}}$ on the B-model according to the prediction of SYZ conjecture. 

\smallskip

 For simplicity, we introduce an extra variable $t\in\R$ and consider the one-parameter family of Laurent polynomials $\{f_t\}$ where $f_t=\sum t^{-\psi(v_i)}z^{v_i}$. Then the large radius limit is achieved by taking $t\rightarrow +\infty$. Moreover, to emphasize the dependence of the form $\omega_c$ on the parameter $t$, we adopt the notation $\omega_{t,c}$ instead. 
 
 \smallskip

\begin{theorem}\label{asymptotics thm}
    The asymptotic behavior of
    \begin{align*}
    	Z^{A,\R_{>0}^d}_c(t^{-\psi(v_1)},...)=\frac{(-1)^d}{(2\pi{\rm i})^{d+1}}\int_{\R_{>0}^d}\omega_{t,c}
    \end{align*}
    when $t\rightarrow +\infty$ is given by
	\begin{align*}
		Z^{A,\R_{>0}^d}_c(t^{-\psi(v_1)},...)=&t^{\psi(c)}\frac{(-1)^{\rk N-\deg c}}{(2\pi{\rm i})^{|\sigma(c)|}|\operatorname{Box}(\sigma(\gamma))|}\cdot\int_{\gamma(c)} t^{\omega}\cdot \hat{\Gamma}_{\gamma(c)} F_{I_c}+o(t^{\psi(c)})
	\end{align*}
	where $\omega=\frac{1}{2\pi {\rm i}}\sum_{i=1}^n\psi(v_i)D_i$, and $\hat{\Gamma}_{\gamma(c)}$ is the Gamma class of $\gamma(c)$ as defined in Definition \ref{def-gamma class}.
\end{theorem}

\smallskip

The proof of this theorem occupies the rest of this section. To begin with, we make a change of coordinates. We denote the moment map on $(\CC^*)^d$ by
\begin{align*}
	\operatorname{Log}_t: (\C^*)^d\rightarrow \R^d, (z_1,\cdots,z_d)\mapsto (\log_t|z_1|,\cdots,\log_t|z_d|)
\end{align*}
and its right-inverse by
\begin{align*}
	i_t: \R^d\rightarrow(\C^*)^d, (y_1,\cdots, y_d)\mapsto (t^{y_1},\cdots, t^{y_d}).
\end{align*}
Note that the positive real locus is identified with $\R^d$ under the map $i_t$. The original integration now becomes 
\begin{align*}
	\int_{\R_{>0}^d}\omega_{t,c}=\int_{\R^d}i_t^*\omega_{t,c}
\end{align*}
with the new coordinates $\{y_i\}$.

\smallskip

We divide the proof into two steps. In the first step (\S\ref{subdivision}), we partition the domain $\R^d$ into smaller sections based on the tropicalization of the Laurent polynomial $f_t$. This allows us to establish a connection between the integration over each section and the volume of specific polytopes. Moving onto the second step (\S\ref{relate to gamma class}) we establish a relationship between these integrals and the integral of Gamma classes on the toric stack $\PP_{\mathbf{\Sigma}}$. A crucial component of the second step is a Duistermaat-Heckman type lemma adapted to our setting that is stated and proved in Appendix \ref{volume}.

\subsection{Subdivision of the domain}\label{subdivision}

For each $i\in\Delta$, we consider the tropicalization $\beta_i$ of the monomial $t^{-\psi(v_i)}z^{v_i}$ defined as
\begin{align*}
	\beta_i:\R^d\rightarrow\R,\quad p\mapsto \langle v_i,p\rangle -\psi(v_i)
\end{align*}
which is an affine function on $\R^d$. Following the idea of \cite{AGIS}, we define\footnote{Our definition differs from the one in \cite{AGIS} by reversing the direction of the inequality, due to the difference between $t\rightarrow +\infty$ and $t\rightarrow 0^+$. }
\begin{align*}
	U^q:=\{p\in\R^d: \beta_i(p)\leq \beta_q(p),\forall i\in\Delta\}
\end{align*}
for any lattice point $q\in\Delta$. Furthermore, for any $q\in\Delta$ and any subset $q\not\in K\subseteq\Delta$ we define a subset $U^{q,K}$ of $U^q$ by
\begin{align*}
	U^{q,K}:=\{p\in\R^d: &\beta_q(p)-\beta_i(p)\in[0,\epsilon],\forall i\in K, \\
	&\beta_q(p)-\beta_i(p)\in[\epsilon,+\infty),\forall i\not\in q\sqcup K\}
\end{align*}
for some fixed small positive number $\epsilon>0$. Intuitively, $U^{q,K}$ is the region where the tropical monomial $\beta_q$ is the largest (hence dominates the asymptotics) and the tropical monomials $\{\beta_k\}_{k\in K}$ are not far behind.

\smallskip

\begin{remark}
	By the standard argument of tropical geometry, we have the following facts about $U^q$ and $U^{q,K}$. Firstly, the set $U^{q,K}$ is non-empty if and only if $q\sqcup K$ forms a cone in $\Sigma$. Additionally, $U^{q,K}$ is unbounded if and only if $q\sqcup K$ is a cone on the boundary, i.e., $\operatorname{relint}(q\sqcup K)\subseteq\partial\Delta$. We will not use the latter fact in rest of the paper so we omit its proof.
\end{remark}

\smallskip

Hence the original integral can be written as a sum
\begin{align*}
	\int_{\R^d}i_t^*\omega_{t,c}=\sum_{q,K}\int_{U^{q,K}}i_t^*\omega_{t,c}
\end{align*}

The first observation is the following lemma which states that only the region $U^{q,K}$ with $\sigma(c)\subseteq q\sqcup K\in \Sigma$ contributes to the leading term when $t\rightarrow+\infty$. Otherwise the growth of the integral over the piece will be $O(t^{\psi(c)-\epsilon})$ for some $\epsilon>0$.

\smallskip

\begin{lemma}
	As $t\rightarrow +\infty$, for $q$ and $K$ with $\sigma(c)\not\subseteq q\sqcup K$ we have $\int_{U^{q,K}}i_t^*\omega_{t,c}=O(t^{\psi(c)-\epsilon})$ for some $\epsilon>0$. If $\sigma(c)\subseteq q\sqcup K\in \Sigma$, then $\int_{U^{q,K}}i_t^*\omega_{t,c}$ is
	\begin{align*}
		(-1)^{\deg c-1}&(\deg c-1)!t^{\psi(c)}(\log t)^d\int_{U^{q,K}}\frac{\prod_{i\in K}(t^{\beta_i-\beta_q})^{c_i}}{\left(1+\sum_{i\in K}t^{\beta_i-\beta_q}\right)^{\deg c}}\prod_{i}dy_i \\
		&+ O(t^{\psi(c)-\epsilon}(\log t)^d).
	\end{align*}
\end{lemma}
\begin{proof}
First we suppose $\sigma(c)\not\subseteq q\sqcup K$. We have
	\begin{align*}
	z^{\bar{c}}&=z^{\sum_{i\in \sigma(c)}c_i \bar{v_i}}=z^{c_q \bar{v_q}}\cdot\prod_{i\in \sigma(c)\cup K\backslash q}z^{c_i \bar{v_i}}\\
&=t^{\sum_{i\in \sigma(c)}c_i\psi(v_i)}(t^{-\psi(v_q)}z^{\bar{v_q}})^{c_q}\cdot\prod_{i\in \sigma(c)\cup K\backslash q}(t^{-\psi(v_i)}z^{\bar{v_i}})^{c_i}\\
&=t^{\sum_{i\in \sigma(c)}c_i\psi(v_i)}(t^{-\psi(v_q)}z^{\bar{v_q}})^{\sum_{i\in q\sqcup K\cup\sigma(c)}c_i}\cdot\prod_{i\in \sigma(c)\cup K\backslash q}(t^{-\psi(v_i)+\psi(v_q)}z^{\bar{v_i}-\bar{v_q}})^{c_i}\\
&=t^{\psi(c)}(t^{-\psi(v_q)}z^{\bar{v_q}})^{\deg c}\cdot\prod_{i\in \sigma(c)\cup K\backslash q}(t^{-\psi(v_i)+\psi(v_q)}z^{\bar{v_i}-\bar{v_q}})^{c_i}
\end{align*}
Similarly we can compute
\begin{align*}
	f_t(z)^{\deg c}&=(\sum_{i\in\Delta}t^{-\psi(v_i)}z^{\bar{v_i}})^{\deg c}\\
	&=(t^{-\psi(v_q)}z^{\bar{v_q}})^{\deg c}\cdot\Bigg(1+\sum_{i\in K}t^{-\psi(v_i)+\psi(v_q)}z^{\bar{v_i}-\bar{v_q}} \\ 
	&\qquad\qquad\qquad\qquad\qquad\qquad\qquad +\sum_{j\not\in q\sqcup K}t^{-\psi(v_j)+\psi(v_q)}z^{\bar{v_j}-\bar{v_q}}\Bigg)^{\deg c}\\
	&=(t^{-\psi(v_q)}z^{\bar{v_q}})^{\deg c}\cdot\left(1+\sum_{i\in K}t^{-\psi(v_i)+\psi(v_q)}z^{\bar{v_i}-\bar{v_q}}+O(t^{-\epsilon})\right)^{\deg c}
\end{align*}
Hence the form $\omega_{t,c}$ is
\begin{align*}
	\omega_{t,c}&=(-1)^{\deg c-1}(\deg c-1)!t^{\psi(c)} \\ 
	& \qquad\qquad\qquad \cdot\frac{\prod_{i\in \sigma(c)\cup K\backslash q}(t^{-\psi(v_i)+\psi(v_q)}z^{\bar{v_i}-\bar{v_q}})^{c_i}}{\left(1+\sum_{i\in K}t^{-\psi(v_i)+\psi(v_q)}z^{\bar{v_i}-\bar{v_q}}+O(t^{-\epsilon})\right)^{\deg c}}\cdot\prod_{i}\frac{dz_i}{z_i}\\
	&=(-1)^{\deg c-1}(\deg c-1)!t^{\psi(c)}\\
	&\qquad\qquad\qquad \cdot\left(\frac{\prod_{i\in \sigma(c)\cup K\backslash q}(t^{-\psi(v_i)+\psi(v_q)}z^{\bar{v_i}-\bar{v_q}})^{c_i}}{\left(1+\sum_{i\in K}t^{-\psi(v_i)+\psi(v_q)}z^{\bar{v_i}-\bar{v_q}}\right)^{\deg c}}+O(t^{-\epsilon})\right)\cdot\prod_{i}\frac{dz_i}{z_i}
\end{align*}
Therefore the pullback $i_t^*\omega_{t,c}$ is
\begin{align*}
	i_t^*\omega_{t,c}&=(-1)^{\deg c-1}(\deg c-1)!t^{\psi(c)}\\
	&\qquad\qquad\qquad \cdot\left(\frac{\prod_{i\in \sigma(c)\cup K\backslash q}(t^{\beta_i-\beta_q})^{c_i}}{\left(1+\sum_{i\in K}t^{\beta_i-\beta_q}\right)^{\deg c}}+O(t^{-\epsilon})\right)\cdot(\log t)^d\prod_{i}dy_i
\end{align*}
So the integration $\int_{U^{q,K}}i_t^*\omega_{t,c}$ is equal to
\begin{align*}
	(-1)^{\deg c-1}&(\deg c-1)!t^{\psi(c)}(\log t)^d\int_{U^{q,K}}\frac{\prod_{i\in \sigma(c)\cup K\backslash q}(t^{\beta_i-\beta_q})^{c_i}}{\left(1+\sum_{i\in K}t^{\beta_i-\beta_q}\right)^{\deg c}}\prod_{i}dy_i \\
		&+ O(t^{\psi(c)-\epsilon}(\log t)^d).
\end{align*}
Notice that if $\sigma(c)\not\subseteq q\sqcup K$, then there exists $i\in\sigma(c)$ such that $\beta_q-\beta_i\geq \epsilon$, where $\epsilon>0$ is the constant used in the definition of $U^{q,K}$, i.e., the nominator of the integrand will contribute a factor of $t^{-\epsilon}$. Therefore the first term in the above expression is also $O(t^{\psi(c)-\epsilon}(\log t)^d)$. By changing $\epsilon$ to a smaller positive number we can adsorb the logarithmic term and get $O(t^{\psi(c)-\epsilon})$.

\smallskip

Now suppose $\sigma(c)\subseteq q\sqcup K$. Then the same computation as above shows that the integral $\int_{U^{q,K}}i_t^*\omega_{t,c}$ is equal to
\begin{align*}
	(-1)^{\deg c-1}&(\deg c-1)!t^{\psi(c)}(\log t)^d\int_{U^{q,K}}\frac{\prod_{i\in K}(t^{\beta_i-\beta_q})^{c_i}}{\left(1+\sum_{i\in K}t^{\beta_i-\beta_q}\right)^{\deg c}}\prod_{i}dy_i \\
		&+ O(t^{\psi(c)-\epsilon}(\log t)^d).
\end{align*}
\end{proof}

\smallskip

According to this lemma we can disregard integrals over $U^{q,K}$ with $\sigma(c)\not\subseteq q\sqcup K$ when computing the leading term of the asymptotic behavior.

We consider the integral
\begin{align}\label{eq1}
	\int_{U^{q,K}}\frac{\prod_{i\in K}(t^{\beta_i-\beta_q})^{c_i}}{\left(1+\sum_{i\in K}t^{\beta_i-\beta_q}\right)^{\deg c}}\prod_{i}dy_i
\end{align}
where $(q,K)$ is a fixed pair with $\sigma(c)\subseteq q\sqcup K$. To simplify the expression, we introduce a change of coordinate on the region $U^{q,K}$. Let us define
\begin{align*}
	b_i:=\beta_q-\beta_i, \text{ for all }i\in K
\end{align*}
and complete $\{b_i\}_{i\in K}$ into an affine coordinate system on $\R^d$ by adding additional covectors $\{e_j\}$. We can then express the standard affine volume form on $\R^d$ in terms of this new system of coordinates:
\begin{align*}
	\prod_i dy_i= r_{q,K}\cdot\prod_{i\in K}db_i\cdot\prod_{j}de_j .
\end{align*}
Thus, the original integral becomes
\begin{align*}
	\int_{U^{q,K}}\frac{t^{-\sum_{i\in K}c_ib_i}}{\left(1+\sum_{i\in K}t^{-b_i}\right)^{\deg c}}\cdot r_{q,K}\cdot\prod_{i\in K}db_i\cdot\prod_{j}de_j .
\end{align*}
Recall that the region $U^{q,K}$ is defined as 
\begin{align*}
	U^{q,K}=\left\{p\in\R^d:b_i\in[0,\epsilon],\forall i\in K;\ \beta_q-\beta_i\in[\epsilon,\infty],\forall i\not\in q\sqcup K\right\}
\end{align*}
We consider the projection $\pi_b: U^{q,K}\rightarrow [0,\epsilon]^{K}$ onto the $(b_i)_{i\in K}$-coordinate plane, and denote the fiber of a fixed $(b_i)_{i\in K}\in[0,\epsilon]^{K}$ by $F^{q,K}((b_i)_{i\in K})$. Then the integral above can be written as an iterated integral
\begin{align*}
	&\int_{U^{q,K}}\frac{t^{-\sum_{i\in K}c_ib_i}}{\left(1+\sum_{i\in K}t^{-b_i}\right)^{\deg c}}\cdot r_{q,K}\cdot\prod_{i\in K}db_i\cdot\prod_{j}de_j\\
	=&\int_{[0,\epsilon]^K}\left(\int_{F^{q,K}((b_i)_{i\in K})}\frac{t^{-\sum_{i\in K}c_ib_i}}{\left(1+\sum_{i\in K}t^{-b_i}\right)^{\deg c}}\cdot r_{q,K}\cdot\prod_{j}de_j \right)\prod_{i\in K}db_i\\
	=&\int_{[0,\epsilon]^K}\left(\frac{t^{-\sum_{i\in K}c_ib_i}}{\left(1+\sum_{i\in K}t^{-b_i}\right)^{\deg c}}\cdot\operatorname{vol}\left(F^{q,K}((b_i)_{i\in K})\right) \right)\prod_{i\in K}db_i .
\end{align*}
By definition, the fiber $F^{q,K}((b_i)_{i\in K})$ is given by
\begin{align*}
	F^{q,K}((b_i)_{i\in K})=\left\{p\in\R^d:\beta_q-\beta_i=b_i,\forall i\in K;\ \beta_q-\beta_i\in[\epsilon,\infty],\forall i\not\in q\sqcup K\right\} .
\end{align*}
To remove the $\epsilon$-dependence, we introduce a new polytope
\begin{align*}
	E^{q,K}((b_i)_{i\in K})=\left\{p\in\R^d:\beta_q-\beta_i=b_i,\forall i\in K;\ \beta_q-\beta_i\in[0,\infty],\forall i\not\in q\sqcup K\right\} .
\end{align*}
and express the fiber $F^{q,K}((b_i)_{i\in K})$ in terms of these new polytopes:
\begin{align*}
	F^{q,K}((b_i)_{i\in K})=E^{q,K}((b_i)_{i\in K})\backslash \left(\bigcup_{j\not\in q\sqcup K}\bigcup_{b_j\in[0,\epsilon]}E^{q,K\sqcup\{j\}}((b_i)_{i\in K},b_j)\right) .
\end{align*}
By inclusion-exclusion principle, we have
\begin{align*}
	\operatorname{vol}(F^{q,K}((b_i)_{i\in K}))=\sum_{J:J\supseteq K,q\not\in J}&(-1)^{|J\backslash K|}\\
	&\cdot\int_{[0,\epsilon]^{J\backslash K}}\operatorname{vol}(E^{q, J}((b_i)_{i\in K},(b_{j}^{\prime})_{j\in J\backslash K}))db^{\prime} .
\end{align*}
Combining these results, the integral \eqref{eq1} becomes
\begin{align*}
	\int_{[0,\epsilon]^K}&\frac{t^{-\sum_{i\in K}c_ib_i}}{\left(1+\sum_{i\in K}t^{-b_i}\right)^{\deg c}}\cdot\\
	&\left(\sum_{J:J\supseteq K,q\not\in J}(-1)^{|J\backslash K|}\int_{[0,\epsilon]^{J\backslash K}}\operatorname{vol}(E^{q, J}((b_i)_{i\in K},(b_{j}^{\prime})_{j\in J\backslash K}))db^{\prime}\right)db .
\end{align*}
By allowing $q$ and $K$ to vary, we obtain
\begin{equation}\label{eq2}
	\begin{aligned}
	\sum_{(J,K,q):J\supseteq K,q\not\in J}(-1)^{|J\backslash K|}\int_{[0,\epsilon]^J}&\frac{t^{-\sum_{i\in K}c_ib_i}}{\left(1+\sum_{i\in K}t^{-b_i}\right)^{\deg c}}\cdot \\ 
	&\operatorname{vol}(E^{q, J}((b_i)_{i\in K},(b_{j}^{\prime})_{j\in J\backslash K}))db^{\prime}db ,
\end{aligned}
\end{equation}
where the sum is taken over all triples $(J,K,q)$ such that $K\subseteq J$ and $q\not\in J$. Note that the summand corresponding to $J$ is zero unless $q\sqcup J$ is a cone in $\Sigma$.

\subsection{Connection to Gamma classes}\label{relate to gamma class}

The goal of this subsection is to reveal the relationship between \eqref{eq2} with the Gamma classes $\Gamma_{\gamma}$ of twisted sectors of the toric stack $\PP_{\mathbf{\Sigma}}$. We adopt a similar approach as presented in \cite{AGIS}.
 
 \smallskip
 
To begin, we utilize the volume formula Proposition \ref{volume and cohomology} to rewrite the sum obtained in the previous subsection as
\begin{equation}\label{eq3}
	\begin{aligned}
	\sum_{(J,K,q):J\supseteq K,q\not\in J}&(-1)^{|J\backslash K|}\int_{[0,\epsilon]^J}\frac{t^{-\sum_{i\in K}c_ib_i}}{\left(1+\sum_{i\in K}t^{-b_i}\right)^{\deg c}}\cdot \\
	&\left(\frac{1}{|\operatorname{Box}(\sigma(\gamma))|}\int_{\gamma}e^{D-\sum_{j\in J}b_j D_j}\cdot\frac{D_{q\sqcup J}}{D_{\sigma(d)}}\cdot F_{I_c}\right)db^{\prime}db
\end{aligned}
\end{equation}
where $D=\sum\psi(v_i)D_i$ and $\gamma:=\gamma(c)$ is the unique twisted sector corresponding to the lattice point $c\in C^{\circ}$. Now we consider the following cohomology class in $H_{\gamma}^*$ obtained by scaling all classes $D_i$ by a factor of $\frac{\log t}{2\pi{\rm i}}$:
\begin{align*}
	P_t=&\sum_{(J,K,q):J\supseteq K,q\not\in J}(-1)^{|J\backslash K|}\left(\frac{\log t}{2\pi{\rm i}}\right)^{|J|+1-|\sigma(c)|}\frac{D_{q\sqcup J}}{D_{\sigma(c)}}\\
	&\cdot\int_{[0,\epsilon]^J}\left(\frac{t^{-\sum_{i\in K}c_ib_i}}{\left(1+\sum_{i\in K}t^{-b_i}\right)^{\deg c}}e^{(\log t)\frac{D}{2\pi{\rm i}}-\sum_{j\in J}b_j (\log t)\frac{D_j}{2\pi{\rm i}}}\right)db^{\prime}db
\end{align*}
Since the integral over $\gamma$ is only relevant to the $\deg=\dim\gamma=d+1-|\sigma(\gamma)|$ part, the expression \eqref{eq3} is equal to
\begin{align}\label{eq4}
	\frac{1}{|\operatorname{Box}(\sigma(\gamma))|}\left(\frac{\log t}{2\pi{\rm i}}\right)^{-(d+1-|\sigma(c)|)}\int_{\gamma} P_t\cdot F_{I_c}
\end{align}
Note that $\deg F_{I_c}=|I_c|$ and $|\sigma(\gamma)|+|I_c|=|\sigma(c)|$.

\smallskip

The goal of the remaining part of this subsection is to prove the following result which relates the cohomology class $P_t$ to the Gamma class $\hat{\Gamma}_{\gamma}$ of the twisted sector $\gamma$.
\begin{proposition}\label{equality of cohomology classes}
The asymptotics of the class $P_t$ is given by
\begin{align*}
	P_t=\frac{(\log t)^{1-|\sigma(c)|}}{(\deg c-1)!}t^{\omega}\widehat{\Gamma}_{\gamma}+O(t^{-\epsilon})
\end{align*}
	where $\omega=\frac{1}{2\pi{\rm i}}\sum\psi(v_i)D_i$.
\end{proposition}

\smallskip

To proceed, we consider the following analytic function in $D_1,\cdots, D_n$, where we think of the variables $D_i$'s as usual complex numbers:
\begin{align*}
	Q_t(D_1,\cdots,&D_n)=\sum_{(J,K,q):J\supseteq K,q\not\in J}(-1)^{|J\backslash K|}\cdot\left(\frac{\log t}{2\pi{\rm i}}\right)^{|J|+1-|\sigma(c)|}e^{(\frac{\log t}{2\pi {\rm i}}) D}\frac{D_{q\sqcup J}}{D_{\sigma(c)}}\\
	&\cdot\int_{[0,\epsilon]^J}\left(\frac{t^{-\sum_{i\in K}c_ib_i}}{\left(1+\sum_{i\in K}t^{-b_i}\right)^{\sum_{i=1}^n(\frac{D_i}{2\pi{\rm i}}+c_i)}}e^{-\sum_{j\in J}b_j (\log t)\frac{D_j}{2\pi{\rm i}}}\right)db^{\prime}db
\end{align*}
The next proposition establishes the relationship between the function $Q$ and the Gamma function $\Gamma$.

\smallskip

\begin{proposition}\label{equality of functions}
	As functions in variables $D_i$'s we have the following identity:
	\begin{align*}
		Q(D_1,\cdots,D_n)=(\log t)^{1-|\sigma(c)|}&(2\pi {\rm i})^{|\sigma(c)|}e^{(\frac{\log t}{2\pi {\rm i}}) D}\\
		&\cdot\frac{\prod_{i=1}^n \frac{D_i}{2\pi{\rm i}}}{D_{\sigma(c)}}\frac{\prod_{i=1}^n\Gamma(\frac{D_i}{2\pi{\rm i}}+c_i)}{\Gamma(\sum_{i=1}^n(\frac{D_i}{2\pi{\rm i}}+c_i))}+O(t^{-\epsilon})
	\end{align*}
\end{proposition}
\begin{proof}

First, we consider a single integral in the definition of $Q$:
\begin{align*}
	\int_{[0,\epsilon]^J}\left(\frac{t^{-\sum_{i\in K}c_ib_i}}{\left(1+\sum_{i\in K}t^{-b_i}\right)^{\sum_{i=1}^n(\frac{D_i}{2\pi{\rm i}}+c_i)}}e^{-\sum_{j\in J}b_j (\log t)\frac{D_j}{2\pi{\rm i}}}\right)db^{\prime}db
\end{align*}
we introduce the change of variables $s_i:=(\log t)b_i$ to rewrite it as
\begin{align*}
	\int_{[0,\epsilon\log t]^J}\left(\frac{e^{-\sum_{i\in K}c_is_i}}{\left(1+\sum_{i\in K}e^{-s_i}\right)^{\sum_{i=1}^n(\frac{D_i}{2\pi{\rm i}}+c_i)}}e^{-\sum_{j\in J}s_j \frac{D_j}{2\pi{\rm i}}}\right)\frac{ds}{(\log t)^{|J|}}
\end{align*}
We now claim that we could replace the region $[0,\epsilon\log t]^J$ of the integral by $[0,\infty)^J$ without changing the leading term of the asymptotics. In other words, we have
	\begin{align*}
	&\int_{[0,\epsilon\log t]^J}\left(\frac{e^{-\sum_{i\in K}c_is_i}}{\left(1+\sum_{i\in K}e^{-s_i}\right)^{\sum_{i=1}^n(\frac{D_i}{2\pi{\rm i}}+c_i)}}e^{-\sum_{j\in J}s_j \frac{D_j}{2\pi{\rm i}}}\right)ds \\
	&=\int_{[0,\infty)^J}\left(\frac{e^{-\sum_{i\in K}c_is_i}}{\left(1+\sum_{i\in K}e^{-s_i}\right)^{\sum_{i=1}^n(\frac{D_i}{2\pi{\rm i}}+c_i)}}e^{-\sum_{j\in J}s_j \frac{D_j}{2\pi{\rm i}}}\right)ds + O(t^{-\epsilon}).
\end{align*}
	To see this, it suffices to observe that the integrand is controlled by 
	\begin{align*}
		e^{-\sum_{i\in K}c_is_i}\cdot e^{-\sum_{j\in J}s_j \frac{D_j}{2\pi{\rm i}}} 
	\end{align*}
	then the claim follows from the fact that $\int_{\epsilon\log t}^{+\infty}e^{-s}ds = O(t^{-\epsilon})$.
	
	\smallskip

Thus it suffices to look at 
\begin{align*}
	\sum_{(J,K,q):J\supseteq K,q\not\in J}(-1)^{|J\backslash K|}\left(\frac{\log t}{2\pi{\rm i}}\right)^{|J|+1-|\sigma(c)|}e^{(\frac{\log t}{2\pi {\rm i}}) D}\frac{D_{q\sqcup J}}{D_{\sigma(c)}}\\
	\cdot\int_{[0,\infty)^J}\left(\frac{e^{-\sum_{i\in K}c_is_i}}{\left(1+\sum_{i\in K}e^{-s_i}\right)^{\sum_{i=1}^n(\frac{D_i}{2\pi{\rm i}}+c_i)}}e^{-\sum_{j\in J}s_j \frac{D_j}{2\pi{\rm i}}}\right)\frac{ds}{(\log t)^{|J|}}
\end{align*}
Note that $\int_{0}^{\infty}e^{-s_j \frac{D_j}{2\pi{\rm i}}}ds_j=(2\pi{\rm i})/D_j$, integrating for all $j\in J\backslash K$, the integral becomes
\begin{align*}
	&\sum_{(J,K,q):J\supseteq K,q\not\in J}(-1)^{|J\backslash K|}\left(\frac{\log t}{2\pi{\rm i}}\right)^{|J|+1-|\sigma(c)|}e^{(\frac{\log t}{2\pi {\rm i}}) D}\frac{D_{q\sqcup J}}{D_{\sigma(c)}}\\
	&\cdot\int_{[0,\infty)^K}\left(\frac{e^{-\sum_{i\in K}c_is_i}}{\left(1+\sum_{i\in K}e^{-s_i}\right)^{\sum_{i=1}^n(\frac{D_i}{2\pi{\rm i}}+c_i)}}e^{-\sum_{i\in K}s_i \frac{D_i}{2\pi{\rm i}}}\prod_{j\in J\backslash K}\frac{2\pi {\rm i}}{D_j}\right)\frac{ds}{(\log t)^{|J|}}
\end{align*}
which is equal to
\begin{align*}
	\sum_{(J,K,q):J\supseteq K,q\not\in J}&(-1)^{|J\backslash K|}\left(\frac{\log t}{2\pi{\rm i}}\right)^{1-|\sigma(c)|}\left(\frac{1}{2\pi {\rm i}}\right)^{|K|}e^{(\frac{\log t}{2\pi {\rm i}}) D}\frac{D_{q\sqcup K}}{D_{\sigma(c)}}\\
	&\cdot\int_{[0,\infty)^K}\left(\frac{e^{-\sum_{i\in K}c_is_i}}{\left(1+\sum_{i\in K}e^{-s_i}\right)^{\sum_{i=1}^n(\frac{D_i}{2\pi{\rm i}}+c_i)}}e^{-\sum_{i\in K}s_i \frac{D_i}{2\pi{\rm i}}}\right)ds
\end{align*}
We rewrite the sum as
\begin{align*}
	\sum_{(K,q):q\not\in K}&\left(\sum_{J:J\supseteq K,q\not\in J}(-1)^{|J\backslash K|}\right)\left(\frac{\log t}{2\pi{\rm i}}\right)^{1-|\sigma(c)|}\left(\frac{1}{2\pi {\rm i}}\right)^{|K|}e^{(\frac{\log t}{2\pi {\rm i}}) D}\frac{D_{q\sqcup K}}{D_{\sigma(c)}}\\
	&\cdot\int_{[0,\infty)^K}\left(\frac{e^{-\sum_{i\in K}c_is_i}}{\left(1+\sum_{i\in K}e^{-s_i}\right)^{\sum_{i=1}^n(\frac{D_i}{2\pi{\rm i}}+c_i)}}e^{-\sum_{i\in K}s_i \frac{D_i}{2\pi{\rm i}}}\right)ds
\end{align*}

\smallskip

Note that for a fixed pair $(K,q)$, the sum
\begin{align*}
	\sum_{J:q\not\in J,K\subseteq J}(-1)^{|J\backslash K|}=\sum_{P\subseteq\{1,\cdots, N\}\backslash(q\sqcup K)}(-1)^{|P|}
\end{align*}
is equal to $(1+(-1))^{|\{1,\cdots, n\}\backslash(q\sqcup K)|}=0$ if $\{1,\cdots, n\}\not=q\sqcup K$ and 1 otherwise. Therefore, the only nonzero term in the sum above corresponds to $K=\{1,\cdots,n\}\backslash q$. Consequently, the sum above can be expressed as
\begin{align*}
	\sum_{q=1}^n&\left(\frac{\log t}{2\pi{\rm i}}\right)^{1-|\sigma(c)|}\left(\frac{1}{2\pi {\rm i}}\right)^{n-1}e^{(\frac{\log t}{2\pi {\rm i}}) D}\frac{\prod_{i=1}^n D_i}{D_{\sigma(c)}}\\
	&\cdot\int_{[0,\infty)^K}\left(\frac{e^{-\sum_{i\in K}c_is_i}}{\left(1+\sum_{i\in K}e^{-s_i}\right)^{\sum_{i=1}^n(\frac{D_i}{2\pi{\rm i}}+c_i)}}e^{-\sum_{i\in K}s_i \frac{D_i}{2\pi{\rm i}}}\right)ds .
\end{align*}

\smallskip

Simplifying further, we obtain
\begin{align*}
	(\log t)^{1-|\sigma(c)|}&(2\pi {\rm i})^{|\sigma(c)|}\sum_{q=1}^n e^{(\frac{\log t}{2\pi {\rm i}}) D}\frac{\prod_{i=1}^n\frac{D_i}{2\pi{\rm i}}}{D_{\sigma(c)}}\\
	&\cdot\int_{[0,\infty)^K}\left(\frac{e^{-\sum_{i\in K}c_is_i}}{\left(1+\sum_{i\in K}e^{-s_i}\right)^{\sum_{i=1}^n(\frac{D_i}{2\pi{\rm i}}+c_i)}}e^{-\sum_{i\in K}s_i \frac{D_i}{2\pi{\rm i}}}\right)ds .
\end{align*}

\smallskip

We substitute $t_i=e^{-s_i}$, leading to the following expression
\begin{align*}
	(\log t)^{1-|\sigma(c)|}&(2\pi {\rm i})^{|\sigma(c)|}e^{(\frac{\log t}{2\pi {\rm i}}) D} \frac{\prod_{i=1}^n \frac{D_i}{2\pi {\rm i}}}{D_{\sigma(c)}}\\
	&\sum_{q\in\{1,\cdots, n\}}\int_{[0,1]^{n-1}}\frac{\prod_{i\not=q}t_i^{\frac{D_i}{2\pi{\rm i}}+c_i}}{\left(1+\sum_{i\not=q}t_i\right)^{\sum_{i=1}^n(\frac{D_i}{2\pi{\rm i}}+c_i)}}\frac{dt}{t}
\end{align*}

\smallskip

Finally we apply Lemma \ref{identity of beta function} to rewrite this expression as
\begin{align*}
	(\log t)^{1-|\sigma(c)|}(2\pi {\rm i})^{|\sigma(c)|}e^{(\frac{\log t}{2\pi {\rm i}}) D}\frac{\prod_{i=1}^n \frac{D_i}{2\pi{\rm i}}}{D_{\sigma(c)}}\cdot\frac{\prod_{i=1}^n\Gamma(\frac{D_i}{2\pi{\rm i}}+c_i)}{\Gamma(\sum_{i=1}^n(\frac{D_i}{2\pi{\rm i}}+c_i))} .
\end{align*}
This concludes the proof of the proposition.
\end{proof}

We have established the desired relationship between the function $Q_t$ and the Gamma functions. Applying them to the cohomology classes $D_i$ along with the fact that $\frac{\prod_{i=1}^n \frac{D_i}{2\pi{\rm i}}}{D_{\sigma(c)}}\cdot\frac{\prod_{i=1}^n\Gamma(\frac{D_i}{2\pi{\rm i}}+c_i)}{\Gamma(\sum_{i=1}^n(\frac{D_i}{2\pi{\rm i}}+c_i))}$ can be written as
\begin{align*}
	\frac{1}{(2\pi {\rm i})^{|\sigma(c)|}}&\frac{\prod_{i\not\in\sigma(c)}\frac{D_i}{2\pi{\rm i}}}{(\deg c-1)!}\prod_{i\not\in\sigma(c)}\Gamma(\frac{D_i}{2\pi{\rm i}})\prod_{i\in\sigma(c)}\Gamma(\frac{D_i}{2\pi{\rm i}}+c_i)\\
	&=\frac{1}{(2\pi {\rm i})^{|\sigma(c)|}(\deg c-1)!}\prod_{i\not\in\sigma(\gamma)}\Gamma(1+\frac{D_i}{2\pi{\rm i}})\prod_{i\in\sigma(\gamma)}\Gamma(\frac{D_i}{2\pi{\rm i}}+\gamma_i)\\
	&=\frac{1}{(2\pi {\rm i})^{|\sigma(c)|}(\deg c-1)!}\cdot\widehat{\Gamma}_{\gamma}
\end{align*}

Hence, the leading term of $\frac{(-1)^{d}}{(2\pi{\rm i})^{d+1}}\int_{\mathbb{R}_{>0}^d}\omega_{t,c}$ is given by (noting that $d+1-\deg c=\rk N-\deg c$)
\begin{align*}
		t^{\psi(c)} \frac{(-1)^{\rk N-\deg c}}{(2\pi {\rm i})^{|\sigma(c)|}|\operatorname{Box}(\sigma(\gamma))|}\int_{\gamma}t^{\omega}\widehat{\Gamma}_{\gamma}\cdot F_{I_c}
\end{align*}
This concludes the proof of Theorem \ref{asymptotics thm}.

\section{Matching the A- and B-brane integral structures}\label{equality}

The goal of this section is to establish the equality between A-brane and B-brane central charges, therefore matching the integral structures on the A- and B-sides defined in Section \ref{central charges}. This is accomplished by utilizing the hypergeometric duality \cite{BHanDuality} as a key ingredient. More precisely, an explicit formula is provided for the pairing between the solution spaces to $\operatorname{bbGKZ}(C,0)$ and its dual $\operatorname{bbGKZ}(C^{\circ},0)$.


\smallskip

\begin{definition}
	For any pair of solutions $(\Phi_c)$ and $(\Psi_d)$ of the systems $\operatorname{bbGKZ}(C,0)$ and $\operatorname{bbGKZ}(C^{\circ},0)$, we define a pairing 
	\begin{align*}
		\langle-,-\rangle: \operatorname{Sol}(\operatorname{bbGKZ}(C,0))\times \operatorname{Sol}(\operatorname{bbGKZ}(C^{\circ},0))\rightarrow\CC
	\end{align*}
	by the following formula
	\begin{align*}
		\langle \Phi,\Psi\rangle = \sum_{c,d,I}\xi_{c,d,I}\Vol_I\left(\prod_{i\in I}x_i\right)\Phi_c\Psi_d
	\end{align*}
	where the coefficient $\xi_{c,d,I}$ is defined as follows. Fix a choice of a generic vector $v\in C^{\circ}$. For a set $I$ of size $\rk N$ we consider the cone
$\sigma_I=\sum_{i\in I} \R_{\geq 0} v_i$. The coefficients $\xi_{c,d,I}$ for $c+d = v_I$ are defined as
\begin{align*}
    \xi_{c,d,I}=\left\{
    \begin{array}{ll}
         & (-1)^{\operatorname{deg}(c)}, \text{ if } \dim\sigma_{I}=\operatorname{rk}N \text{ and both } c+\varepsilon v\text{ and }d-\varepsilon v\in\sigma_{I}^{\circ} \\
         & 0, \text{ otherwise}.
    \end{array}
    \right.
\end{align*}
Here the condition needs to hold for all sufficiently small positive number $\varepsilon>0$. 
\end{definition}

\smallskip

The main result in \cite{BHanDuality} states that this pairing is non-degenerate.

\smallskip

\begin{theorem}
	For any pair of solutions $(\Phi_c)$ and $(\Psi_d)$ of the systems $\operatorname{bbGKZ}(C,0)$ and $\operatorname{bbGKZ}(C^{\circ},0)$, the pairing $\langle \Phi,\Psi\rangle$ is a constant. Furthermore, the constant pairing $\langle \Gamma,\Gamma^{\circ}\rangle$ of the cohomology-valued Gamma series is equal  to the inverse of the Euler characteristic pairing $\chi:H_{\mathrm{orb}}^*\otimes H_{\mathrm{orb},c}^*\rightarrow\CC$ in the large radius limit. In particular, $\langle-,-\rangle$ is non-degenerate.
\end{theorem}

\smallskip

We now combine the computation in section \ref{asymptotic} together with the hypergeometric duality to obtain the equality between A-brane and B-brane central charges. Specifically, we begin by proving the equality for the case of structure sheaf $\mathcal{O}_{\PP_{\mathbf{\Sigma}}}$ and its mirror cycle $\R_{\geq 0}^d$. 

\smallskip

To start with, we recall the asymptotic behavior of the Gamma series $\Gamma$ that was computed in \cite{BHanDuality}.

\smallskip

\begin{lemma}\label{asymptotics of gamma series}
	Let $t\rightarrow +\infty$, then for lattice point $c\in C$ and $\gamma\not=\gamma^{\vee}(c)$, the summand of $\Gamma_c(t^{-\psi(v_1)}x_1,\cdots,t^{-\psi(v_n)}x_n)$ is $o(t^{\psi(c)})$. For $\gamma=\gamma^{\vee}(c)$, we have
\begin{align*}
\Gamma_c(t^{-\psi(v_1)} x_1, \ldots)
=
t^{\psi(c)} \prod_{i=1}^n {\rm e}^{\frac {D_i}{2\pi{\rm i}}( \log x_i - \psi(v_i) \log t)}
\prod_{i=1}^{n}\frac {x_i^{-c_i}}{\Gamma(1-c_i +\frac {D_i}{2\pi{\rm i}} )}
(1+o(1)).
\end{align*}
\end{lemma}
\begin{proof}
	See \cite[Lemma 3.10]{BHanDuality}.
\end{proof}

\smallskip

\begin{theorem}
	The A-brane central charge associated to the positive real locus $(\R_{\geq 0})^d$ coincides with the B-brane central charge associated to the structure sheaf $\mathcal{O}_{\PP_{\mathbf{\Sigma}}}$.
\end{theorem}
\begin{proof}
Throughout this proof, we will denote an interior lattice point by $d\in C^{\circ}$ and denote the rank of the lattice $N$ by $\rk N$.

	Consider the pairing
	\begin{align*}
		\langle \Gamma,Z^{A,\R_{>0}^d} \rangle=\sum_{c,d,I}\xi_{c,d,I}\Vol_I(\prod_{i\in I}x_i)\Gamma_c\cdot Z^{A,\R_{>0}^d}_d\in H_{\mathrm{orb}}^*(\PP_{\mathbf{\Sigma}})
	\end{align*}
	we look at the component corresponding to a fixed twisted sector $\gamma$. Combining Theorem \ref{asymptotics thm}, Lemma \ref{asymptotics of gamma series}, by an argument similar to the proofs of \cite[Proposition 3.12, 3.13]{BHanDuality}, the asymptotic behavior of
	\begin{align*}
		\prod_{i=1}^n(t^{-\psi(v_i)})\Gamma_{c,\gamma}(t^{-\psi(v_i)})Z^{A,\R_{>0}^d}_d(t^{-\psi(v_i)})
	\end{align*}
	is given by $o(1)$ unless $\gamma=\gamma^{\vee}(c)=\gamma(d)$ and both $I_c$ and $I_d$ are cones in $\Sigma$, in which case the asymptotic behavior is
	\begin{align*}
		o(1)+\frac{1}{(2\pi{\rm i})^{|I_c|}}\cdot \frac{D_{I_c}}{\hat{\Gamma}_{\gamma}}\prod_{i=1}^n {\rm e}^{\frac{D_i}{2\pi {\rm i}}(-\psi(v_i)\log t)}\frac{(-1)^{\rk N-\deg d}}{(2\pi{\rm i})^{|\sigma(d)|}|\operatorname{Box}(\sigma(\gamma))|}\int_{\gamma} t^{\omega}\hat{\Gamma}_{\gamma} F_{I_d}.
	\end{align*}
	Since $\langle \Gamma,\Psi \rangle$ is a constant, taking the constant term we get
	\begin{align*}
		\langle \Gamma,\Psi \rangle_{\gamma}&=\sum_{c,d,I}\xi_{c,d,I}\Vol_I \frac{D_{I_c}}{\hat{\Gamma}_{\gamma}}\frac{(-1)^{\rk N-\deg d}}{(2\pi{\rm i})^{|I_c|+|\sigma(d)|}|\operatorname{Box}(\sigma(\gamma))|}\int_{\gamma}\hat{\Gamma}_{\gamma} F_{I_d}\\
		&=\sum_{c,d,I}\xi_{c,d,I}\Vol_I \frac{D_{I_c}}{\hat{\Gamma}_{\gamma}}\frac{(-1)^{\rk N-\deg d}}{(2\pi{\rm i})^{\rk N}|\operatorname{Box}(\sigma(\gamma))|}\int_{\gamma^{\vee}} (-1)^{\dim\gamma^{\vee}-|I_d|} (\hat{\Gamma}_{\gamma})^* F_{I_d}\\
		&=\sum_{c,d,I}\xi_{c,d,I}\Vol_I \frac{D_{I_c}}{\hat{\Gamma}_{\gamma}}\frac{(-1)^{\rk N-\deg d}}{(2\pi{\rm i})^{\rk N}|\operatorname{Box}(\sigma(\gamma))|}\\
		&\quad\quad\quad \cdot\int_{\gamma^{\vee}} (2\pi{\rm i})^{|\sigma(\gamma)|}(-1)^{\deg\gamma^{\vee}+\dim\gamma^{\vee}-|I_d|} \frac{F_{I_d}}{\hat{\Gamma}_{\gamma^{\vee}}}\mathrm{Td}(\gamma^{\vee})\\
		&=\frac{1}{(2\pi{\rm i})^{\rk N}}\sum_{c,d,I}\xi_{c,d,I}\Vol_I (2\pi{\rm i})^{|\sigma(\gamma)|} \frac{D_{I_c}}{\hat{\Gamma}_{\gamma}}\frac{1}{|\operatorname{Box}(\sigma(\gamma))|}\int_{\gamma^{\vee}} \frac{F_{I_d}}{\hat{\Gamma}_{\gamma^{\vee}}}\mathrm{Td}(\gamma^{\vee})\\
		&=\frac{1}{(2\pi{\rm i})^{\rk N}}\sum_{c,d,I}\xi_{c,d,I}\Vol_I (2\pi{\rm i})^{|\sigma(\gamma)|} \frac{D_{I_c}}{\hat{\Gamma}_{\gamma}}\cdot(1,\frac{F_{I_d}}{\hat{\Gamma}_{\gamma^{\vee}}})_{\mathrm{orb},\gamma}
	\end{align*}
	Here we used $\deg\gamma^{\vee}=|\sigma(\gamma)|+|I_d|-\deg d=|\sigma(d)|-\deg d$, therefore
	\begin{align*}
		(-1)^{\rk N-\deg d+\deg\gamma^{\vee}+\dim\gamma^{\vee}-|I_d|}&=(-1)^{\rk N+|\sigma(d)|+\rk N-|\sigma(\gamma)|-|I_d|}\\
		&=(-1)^{\rk N+|\sigma(d)|+\rk N-|\sigma(d)|}=1
	\end{align*}
	By \cite[Theorem 4.2]{BHanDuality} the following class in $H_{\mathrm{orb}}^*\otimes H_{\mathrm{orb},c}^*$
	\begin{align*}
		\frac{1}{(2\pi{\rm i})^{\rk N}}\bigoplus_{\gamma}\sum_{c,d,I}\xi_{c,d,I}\Vol_I (2\pi{\rm i})^{|\sigma(\gamma)|} \frac{D_{I_c}}{\hat{\Gamma}_{\gamma}}\otimes\frac{F_{I_d}}{\hat{\Gamma}_{\gamma^{\vee}}}
	\end{align*}
	is inverse to the Euler characteristic pairing, therefore for any $\gamma$ we have
	\begin{align*}
		\langle \Gamma,Z^{A,\R_{>0}^d} \rangle_{\gamma}=1_{\gamma}
	\end{align*}
	so $\langle \Gamma,Z^{A,\R_{>0}^d} \rangle=\bigoplus_{\gamma}1_{\gamma}=ch(\mathcal{O}_{\PP_{\mathbf{\Sigma}}})$, i.e., $Z^{A,\R_{>0}^d}$ corresponds to the structure sheaf $\mathcal{O}_{\PP_{\mathbf{\Sigma}}}$.
\end{proof}

\smallskip

We have completed the proof of Theorem \ref{main theorem} for the case of structure sheaf $\mathcal{O}_{\PP_{\mathbf{\Sigma}}}$. Next, we consider an arbitrary line bundle $\mathcal{L}=\mathcal{O}(\sum_{i=1}^n a_i D_i)$ corresponding to a torus-invariant divisor $\sum_{i=1}^n a_iD_i$. The mirror cycle $\operatorname{mir}(\mathcal{L})$ of $\mathcal{L}$ is constructed from $\R_{>0}^d$ by monodromy. More precisely, the divisor $\sum_{i=1}^n a_iD_i$ defines a loop in the stringy K\"ahler moduli space of $\PP_{\mathbf{\Sigma}}$ by
\begin{align}\label{loop}
	\phi:[0,1]\rightarrow \CC^n,\quad \theta\mapsto (e^{-2\pi{\rm i}a_1\theta},\cdots,e^{-2\pi{\rm i}a_n\theta})
\end{align}
We denote the Laurent polynomial corresponding to $\phi(\theta)$ by $f^{(\theta)}$, then we have a family of hypersurfaces $Z_{f^{(\theta)}}$ in $(\C^*)^d$, where $Z_{f^{(1)}}=Z_{f^{(0)}}=Z_{f}$. We then define the mirror cycle of $\mathcal{L}$ to be the parallel transport of $\R_{>0}^d$ along this loop.

\smallskip

Note that the set of $\mathrm{mir}(\mathcal{L})$ for all $\mathcal{L}\in K_0(\PP_{\mathbf{\Sigma}})$ generate a sublattice of the integral homology group $H_{d}\left((\CC^*)^{d}\backslash Z_f,\Z\right)$ that has the correct rank $\mathrm{vol}(\Delta)$, which is a direct consequence of the following result.

\smallskip

\begin{corollary}\label{monodromy}
	For any $\mathcal{L}=\mathcal{O}(\sum_{i=1}^n a_i D_i)\in K_0(\PP_{\mathbf{\Sigma}})$, the A- and B-brane central charges coincide:
	\begin{align*}
		Z^{A,\mathrm{mir}(\mathcal{L})}=Z^{B,\mathcal{L}}.
	\end{align*}
	As a consequence, the A- and B-brane integral structures of the Hori-Vafa mirrors, defined by $H_{d}\left((\CC^*)^{d}\backslash Z_f,\Z\right)$ and 
 $K_0(\PP_{\Sigma},\Z)$ respectively, coincide.

\end{corollary}
\begin{proof}
	It suffices to compare the monodromy along the loop \eqref{loop} on both sides. Recall that the Gamma series is given by
	\begin{align*}
\Gamma_c^{\circ}(x_1,\ldots,x_n) =\bigoplus_{\gamma} \sum_{l\in L_{c,\gamma}} \prod_{i=1}^n 
\frac {x_i^{l_i+\frac {D_i}{2\pi \ii }}}{\Gamma(1+l_i+\frac {D_i}{2\pi \ii })}\left(\prod_{i\in\sigma}D_i^{-1}\right)F_{\sigma}
\end{align*}
and the monodromy comes from the term $\prod_{i=1}^n x_i^{l_i+\frac {D_i}{2\pi \ii }}=\prod_{i=1}^n e^{(l_i+\frac {D_i}{2\pi \ii })\log x_i}$.

Fix $c$ and $\gamma$, when $\theta$ goes from 0 to 1, the $x_i$ goes around the origin $a_i$ times clockwisely and hence the original $\log x_i$ now becomes $\log x_i - 2\pi {\rm i}a_i$, therefore contributes an extra factor $e^{-a_i(2\pi {\rm i}l_i+ D_i)}$. Take product over all $i=1,\cdots, n$, this is $e^{-\sum_{i}a_i(2\pi {\rm i}l_i+ D_i)}$. By definition of $l\in L_{c,\gamma}$, we have $l_i\equiv\gamma_i\mod\Z$, therefore the factor is equal to $e^{-\sum_{i}a_i(2\pi {\rm i}\gamma_i+ D_i)}$, which is exactly the Chern character $ch_{\gamma}(\mathcal{O}(-\sum_{i=1}^n a_i D_i))$.

Consequently, the effect of the monodromy on the Gamma series is to multiply it by $ch(\mathcal{O}(-\sum_{i=1}^n a_i D_i))$. Then the B-brane central charge is obtained by composing with $\chi(ch(\mathcal{O}_{\PP_{\mathbf{\Sigma}}}),-)$, which is
\begin{align*}
	\chi\left(1,ch(\mathcal{O}(-\sum_{i=1}^n a_i D_i))\cdot\Gamma_c^{\circ}\right)=\chi\left(ch(\mathcal{O}(\sum_{i=1}^n a_i D_i)),\Gamma_c^{\circ}\right)
\end{align*}
by Proposition \ref{def-eulerpairing}. This is exactly the central charge $Z^{B,\mathcal{O}(\sum_{i=1}^n a_i D_i)}$. It then follows directly from the construction of the mirror cycle of $\mathcal{O}(\sum a_i D_i)$ that the monodromy on the A-side matches with the monodromy on the B-side.
\end{proof}

\appendix

\section{Residual volume and orbifold cohomology}\label{volume}

In this appendix we establish a result that relates the residual volume of the polytopes $E^{q,J}((b_j)_{j\in J})$ (for precise definitions of residual volume and the polytopes, see section \ref{asymptotic}) with certain orbifold cohomology classes with compact support of the toric Deligne-Mumford stack $\PP_{\mathbf{\Sigma}}$. It could be considered as a replacement of the Duistermaat-Heckman lemma used in \cite{AGIS} adapted to our setting. We use the same notations from section \ref{asymptotic} except we denote a twisted sector corresponding to $\gamma\in\operatorname{Box}(\Sigma)$ by $\PP_{\mathbf{\Sigma}/\gamma}$ to avoid potential confusion.

\smallskip

\begin{proposition}\label{volume and cohomology}
	The residual volume $\operatorname{vol}(E^{q,J}((b_j)_{j\in J}))$ is equal to
	\begin{align*}
		\frac{1}{|\operatorname{Box}(\sigma(\gamma))|}\int_{\PP_{\mathbf{\Sigma}/\gamma}}e^{D-\sum_{j\in J}b_jD_j}\frac{D_{q\sqcup J}}{D_{\sigma(c)}}F_{I_c}
	\end{align*}
	where $D:=\sum_{i}\psi(v_i)D_i$, and $\gamma:=\gamma(c)$ is the unique twisted sector corresponding to the interior lattice point $c\in C^{\circ}$.
\end{proposition}

\smallskip

We begin with a review of the well-known results on the relationship between line bundles on toric varieties and the associated polytopes in \S\ref{line bundles and polytopes} and provide the proof of Proposition \ref{volume and cohomology} in \S\ref{proof of the volume formula}.

\subsection{Line bundles on toric varieties and their associated polytopes}\label{line bundles and polytopes}

We briefly review the classical correspondence between line bundles on toric varieties and their associated polytopes following \cite{Cox}.

\smallskip

Again, we denote by $\PP_{\Sigma}$ the toric variety corresponding to a fan $\Sigma$, and $D=\sum_{\rho}a_{\rho}D_{\rho}$ be a Cartier divisor on $\PP_{\Sigma}$, where $D_{\rho}$'s are the torus-invariant divisors, and we denote the primitive generators of the correponding rays in the fan by $v_{\rho}$. The associated polytope $P_D$ of the line bundle $\mathcal{O}_{\PP_{\Sigma}}(D)$ is defined as\footnote{There is a difference of signs in our definition with the one in \cite{Cox}.}
\begin{align*}
	P_D:=\{m\in M_{\R}: \langle m, v_{\rho}\rangle  +a_{\rho}\geq 0,\ \forall\rho\}
\end{align*}
It is a well-known fact that the dimension of the global section of $\mathcal{O}_{\PP_{\Sigma}}(D)$ is equal to the number of lattice points in the polytope $P_D$. In fact, we have
\begin{align*}
	\Gamma(\PP_{\Sigma},\mathcal{O}_{\PP_{\Sigma}}(D))=\bigoplus_{m\in P_D\cap M}\CC\cdot\chi^m
\end{align*}
In this case the polytope $P_D$ is of full-dimension.

\smallskip

This correspondence could be generalized further to the case where the polytope is not of full-dimension. In this case, the corresponding sheaf is not a line bundle on $\PP_{\Sigma}$, but the restriction of a line bundle to a closed subvariety.

\smallskip

Finally, suppose we have a sheaf of the form $\mathcal{O}_{D^{\prime}}(D)$ where $D^{\prime}$ and $D$ are torus-invariant divisors and its associated (non-full-dimensional) polytope $P$. Suppose further that this sheaf is nef. Then by Demazure vanishing theorem (see \cite[Theorem 9.2.3]{Cox}) all higher cohomology of $\mathcal{O}_{D^{\prime}}(D)$ vanishes. Consequently, we have
\begin{align*}
	\chi(\PP_{\Sigma},\mathcal{O}_{D^{\prime}}(D))=\chi(D^{\prime},\mathcal{O}_{D^{\prime}}(D))=\dim H^0(D^{\prime},\mathcal{O}_{D^{\prime}}(D))=|P\cap M|
\end{align*}
i.e., the Euler characteristic of $\mathcal{O}_{D^{\prime}}(D)$ is equal to the number of lattice points in the polytope $P$.

\subsection{Proof of Proposition \ref{volume and cohomology}}\label{proof of the volume formula}
Before we start, we remark here that it suffices to prove the statement for the case where $q\sqcup J$ is a cone in the fan $\Sigma$ because otherwise the polytope $E^{q,J}((b_j)_{j\in J})$ is empty and the right hand side of the equality is zero due to the factor $D_{q\sqcup J}$.

\smallskip

We divide the proof into four steps. 

\smallskip

{\bf Step 1.}
	Recall that the polytope $E^{q,J}((b_j)_{j\in J})$ is defined as
	\begin{align*}
	E^{q,J}((b_i)_{i\in J})=\left\{p\in\R^d:\beta_q-\beta_i=b_i,\forall i\in J;\ \beta_q-\beta_i\in[0,\infty],\forall i\not\in q\sqcup J\right\}
    \end{align*}
    where $\beta_i:\R^d\rightarrow\R$ is a linear function defined as $p\mapsto \langle v_i, p \rangle-\psi(v_i)$. Without loss of generality, we assume that all $b_i$'s are rational numbers. The defining equalities and inequalities of $E^{q,J}((b_j)_{j\in J})$ could be rewritten as
    \begin{align*}
    	\langle v_q-v_i, p\rangle + \psi(v_i)-\psi(v_q)-b_i = 0
    \end{align*}
    for $i\in J$ and
    \begin{align*}
    	\langle v_q-v_i, p\rangle + \psi(v_i)-\psi(v_q) \geq 0
    \end{align*}
    for $i\not\in q\sqcup J$. If we denote
    \begin{align*}
    	D:&=\sum_{i\in J}(\psi(v_i)-\psi(v_q)-b_i)D_i + \sum_{i\not\in q\sqcup J}(\psi(v_i)-\psi(v_q))D_i \\ 
    	&=\sum_{i\in\operatorname{Star}(q)}\psi(v_i)D_i - \sum_{i\in J}b_i D_i
    \end{align*}
    then by the discussion in \S\ref{line bundles and polytopes}, we have
    \begin{align*}
    	\chi(\PP_{\Sigma/q},\mathcal{O}_{D_J}(lD))=|l\cdot E^{q,J}((b_i)_{i\in J})\cap M|
    \end{align*}
    where $l$ is any integer number that makes $l\cdot E^{q,J}((b_i)_{i\in J})$ into a lattice polytope (the existence is due to the rationality of $b_i$'s), and $D_J$ is the closed subvariety of $\PP_{\Sigma/q}$ corresponding to the cone $J$ in the quotient fan $\Sigma/q$. This can be further expressed as
    \begin{align*}
    	|l\cdot E^{q,J}((b_i)_{i\in J})\cap M|=\chi(\PP_{\Sigma/(q\sqcup J)},\mathcal{O}_{\PP_{\Sigma/(q\sqcup J)}}(lD))
    \end{align*}
    
\medskip
    
{\bf Step 2.}
Denote the canonical map from the smooth toric Deligne-Mumford stack $\PP_{\mathbf{\Sigma}}$ to its coarse moduli space (which is a simplicial toric variety) $\PP_{\Sigma}$ by $\pi$. We denote the line bundle on $\PP_{\mathbf{\Sigma}}$ defined by the same support function with $D$ by $\mathcal{O}_{\PP_{\mathbf{\Sigma}}}(D)$. Since in the definition of $\PP_{\mathbf{\Sigma}}$ the additional data of a vector on each ray of the stacky fan is chosen to be the primitive generator, we know that the pushforward of $\mathcal{O}_{\PP_{\mathbf{\Sigma}}}(D)$ is exactly $\mathcal{O}_{\PP_{\Sigma}}(D)$. On the other hand, it is a well-known result (see e.g.,\cite[Definition 4.1, Example 8.1]{stack}) for a tame Deligne-Mumford stack $\mathcal{X}$ with coarse moduli space $X$, the canonical map $\pi:\mathcal{X}\rightarrow X$ is cohomologically affine. This implies that $H^i(\mathcal{X},\mathcal{F})$ is equal to $H^i(X,\pi_*\mathcal{F})$ for all $i>0$ and any coherent sheaf $\mathcal{F}$. Apply this fact to our situation, we get $\chi(\PP_{\mathbf{\Sigma}/(q\sqcup J)},\mathcal{O}_{\PP_{\mathbf{\Sigma}/(q\sqcup J)}}(lD))=\chi(\PP_{\Sigma/(q\sqcup J)},\mathcal{O}_{\PP_{\Sigma/(q\sqcup J)}}(lD))$. Thus, we have
\begin{align*}
	\chi(\PP_{\mathbf{\Sigma}/(q\sqcup J)},\mathcal{O}_{\PP_{\mathbf{\Sigma}/(q\sqcup J)}}(lD))=|l\cdot E^{q,J}((b_i)_{i\in J})\cap M|
\end{align*}

\medskip

{\bf Step 3.}

Then we apply Corollary \ref{formula of euler characteristic}, we obtain
\begin{align*}
	|l\cdot E^{q,J}((b_i)_{i\in J})\cap M| &= \sum_{\gamma\in\operatorname{Box}(\Sigma/(q\sqcup J))}\frac{1}{|\operatorname{Box}(\sigma(\gamma))|}\\
	&\qquad\qquad\qquad \cdot\int_{\gamma}ch_{\gamma}^c(l\cdot \mathcal{O}_{\PP_{\mathbf{\Sigma}/(q\sqcup J)}}(lD))\operatorname{Td}(\gamma)\\
	&=\sum_{\gamma\in\operatorname{Box}(\Sigma/(q\sqcup J))}\frac{1}{|\operatorname{Box}(\sigma(\gamma))|}\int_{\gamma}e^{lD}\cdot \operatorname{Td}(\gamma)
\end{align*}
note that since $q\sqcup J$ is an interior cone (because it contains an interior cone $\sigma(c)$ as a subcone), the quotient fan $\Sigma/(q\sqcup J)$ is complete, hence $\PP_{\mathbf{\Sigma}/(q\sqcup J)}$ is compact therefore $K_0$ and $K_0^c$ (and the corresponding Chern characters) coincide.

\smallskip

The {\it affine volume}\footnote{Note that the affine volume differs with the residual volume $\operatorname{vol}E^{q,J}((b_i)_{i\in J})$ by a factor of the index of $b_i$'s, namely the index of the sublattice spanned by $b_i$'s inside the standard lattice $\Z^{d}$, see step 4.} $\operatorname{vol}_{\mathrm{aff}}E^{q,J}((b_i)_{i\in J})$ is computed by
\begin{align*}
	\operatorname{vol}_{\mathrm{aff}}E^{q,J}((b_i)_{i\in J})&=\lim_{l\rightarrow\infty}\frac{|l\cdot E^{q,J}((b_i)_{i\in J})\cap M|}{l^{\dim E^{q,J}((b_i)_{i\in J})}}\\
	&=\sum_{\gamma\in\operatorname{Box}(\Sigma/(q\sqcup J))}\frac{1}{|\operatorname{Box}(\sigma(\gamma))|}\lim_{l\rightarrow\infty}\int_{\gamma}\frac{e^{lD}}{l^{\rk N-1-|J|}}\cdot \operatorname{Td}(\gamma) .
\end{align*}
In the last step, we used $\dim E^{q,J}((b_i)_{i\in J})=\rk N-1-|J|$. 

\smallskip

Now, we claim that the only nonzero term in this sum is the $\gamma=0$ term. To see this, we expand the $e^{lD}$ and the Todd class $\operatorname{Td}(\gamma)$ as sums of homogeneous components:
\begin{align*}
	\int_{\gamma}\frac{e^{lD}}{l^{\rk N-1-|J|}}\cdot \operatorname{Td}(\gamma) = \int_{\gamma}\sum_{i,j=0}^{\infty}\frac{1}{l^{\rk N-1-|J|}}\frac{l^i D^i}{i!}\operatorname{Td}(\gamma)_j
\end{align*}
where $\operatorname{Td}(\gamma)_j$ denotes the degree $j$ part of the Todd class of $\gamma$. By definition, the only nonzero contribution comes from terms with $(i,j)$ such that
\begin{align*}
	\deg(D^i)+\deg(\operatorname{Td}(\gamma)_j)=i+j
\end{align*}
is exactly equal to $\dim(\gamma)=\rk N-1-|J|-|\sigma(\gamma)|$. On the other hand, if $i<\rk N-1-|J|$, the integral will be killed by taking limit $l\rightarrow\infty$ due to the factor of $\frac{1}{l^{\rk N-1-|J|}}$, thus $i\geq \rk N-1-|J|$. Combining these two observations, we can deduce that in order to have nonzero contribution, we must have
\begin{align*}
	\rk N-1-|\sigma(\gamma)|-|J|=i+j\geq \rk N-1-|J|+j
\end{align*}
which simplifies to $|\sigma(\gamma)|\leq -j$. This forces $j=0$, $i=\rk N-1-|J|$ and $\sigma(\gamma)=\emptyset$, i.e, $\gamma=0$. Hence the claim is proved.
Therefore we have
\begin{align*}
	\operatorname{vol}_{\mathrm{aff}}E^{q,J}((b_i)_{i\in J})=\int_{\PP_{\mathbf{\Sigma}/(q\sqcup J)}}\frac{D^{\rk N-1-|J|}}{(\rk N-1-|J|)!}=\int_{\PP_{\mathbf{\Sigma}/(q\sqcup J)}}e^{D}
\end{align*}

\medskip

{\bf Step 4.}
The affine and residual volume of the polytope $E^{q,J}((b_i)_{i\in J})$ are related by the following equation:
\begin{align*}
	\operatorname{vol}_{\mathrm{aff}}E^{q,J}((b_i)_{i\in J}) = (\text{index of }b_i's)\cdot \operatorname{vol}E^{q,J}((b_i)_{i\in J})
\end{align*}
The index of $b_i$'s is exactly equal to $|\operatorname{Box}(q\sqcup J)|$. Therefore, we have
\begin{align*}
	\operatorname{vol}E^{q,J}((b_i)_{i\in J})&=\frac{1}{|\operatorname{Box}(q\sqcup J)|}\int_{\PP_{\mathbf{\Sigma}/(q\sqcup J)}}e^{D}\\
	&=\int_{\PP_{\mathbf{\Sigma}}}e^{D}\cdot F_{q\sqcup J}\\
	&=\frac{1}{|\operatorname{Box}(\sigma(\gamma))|}\int_{\PP_{\mathbf{\Sigma}/\sigma(\gamma)}}e^{D}\cdot F_{q\sqcup J\backslash\sigma(\gamma)}\\
	&=\frac{1}{|\operatorname{Box}(\sigma(\gamma))|}\int_{\PP_{\mathbf{\Sigma}/\sigma(\gamma)}}e^{D}\cdot \frac{D_{q\sqcup J}}{D_{\sigma(c)}}F_{I_c}
\end{align*}
where the second equality follows from the fact that the ratio between the volume of a cone in the fan $\Sigma$ and that of the corresponding quotient cone in the quotient fan $\Sigma/(q\sqcup J)$ is equal to $|\operatorname{Box}(q\sqcup J)|$. The third equality holds for similar reasons. The last equality is a consequence of the relations in the orbifold cohomology space. This concludes the proof.

\section{An identity of Beta function}\label{beta function}
In this appendix, we prove an identity of multivariate beta functions that is used in the computation in \S\ref{asymptotic}. Recall that the usual beta function $B(a,b)$ is defined as
\begin{align*}
	B(a,b):=\int_{0}^1 t^{a-1}(1-t)^{b-1}dt
\end{align*}
There is a well-known relation between beta function and gamma function $B(a,b)=\Gamma(a)\Gamma(b)/\Gamma(a+b)$. There is an equivalent definition of $B(a,b)$ given by
\begin{align*}
	B(a,b)=\int_{0}^{\infty}\frac{t^{a-1}}{(1+t)^{a+b}}dt
\end{align*}

\smallskip

Similarly, we can define the multivariable beta function by
\begin{align*}
	B(a_1,\cdots,a_n)=\int_{[0,\infty)^{n-1}}\frac{t_1^{a_1}\cdots t_{n-1}^{a_{n-1}}}{(1+t_1+\cdots+t_{n-1})^{a_1+\cdots+a_n}}dt_1\cdots dt_{n-1}
\end{align*}
and there is an identity $B(a_1,\cdots,a_n)=\Gamma(a_1)\cdots\Gamma(a_n)/\Gamma(a_1+\cdots+a_n)$.

\smallskip

\begin{lemma}\label{identity of beta function}
We have the following identity:
	\begin{align*}
		B(a_1,\cdots,a_n)=\sum_{q=1}^n \int_{[0,1]^{n-1}}\frac{\prod_{i\not=q}t_i^{a_i-1}}{\left(1+\sum_{i\not=q}t_i\right)^{\sum_{i=1}^n a_i}}\prod_{i\not=q}dt_i
	\end{align*}
\end{lemma}
\begin{proof}
	We treat the $q=n$ term and the other terms separately. The $q=n$ term could be written as
	\begin{align*}
		\int_{[0,1]^{n-1}}\frac{\prod_{i=1}^{n-1}t_i^{a_i-1}}{\left(1+\sum_{i=1}^{n-1}t_i\right)^{\sum_{i=1}^n a_i}}\prod_{i=1}^{n-1}dt_i = \sum_{q=1}^{n-1} \int_{A_q}\frac{\prod_{i=1}^{n-1}t_i^{a_i-1}}{\left(1+\sum_{i=1}^{n-1}t_i\right)^{\sum_{i=1}^n a_i}}\prod_{i=1}^{n-1}dt_i
	\end{align*}
	where $A_q$ is the region defined by $0\leq t_q\leq 1$ and $0\leq t_j\leq t_q$ for $j\not=q$.

	\smallskip
	
	 When $q\not=n$, we introduce the change of coordinate given by $t_i\rightarrow \frac{t_i}{t_n}$ for $i\not=n$ and $t_n\rightarrow \frac{1}{t_n}$, then rename the variable $t_n$ to $t_q$. An elementary computation shows that the integral becomes
	\begin{align*}
		\int_{B_q}\frac{\prod_{i=1}^{n-1}t_i^{a_i-1}}{\left(1+\sum_{i=1}^{n-1}t_i\right)^{\sum_{i=1}^n a_i}}\prod_{i=1}^{n-1}dt_i
	\end{align*}
	where $B_q$ is the region defined by $t_q\geq 1$ and $0\leq t_j\leq t_q$ for $j\not=q$. Therefore the original sum can be written as the integral over the union $\bigcup_{q=1}^{n-1}(A_q\cup B_q)$. Now the result follows from the observation that $A_q\cup B_q$ is the region defined by $t_q\geq 0$ and $0\leq t_j\leq t_q$ for $j\not=q$, and $\bigcup_{q=1}^{n-1}(A_q\cup B_q)$ is exactly $[0,\infty)^{n-1}$.
\end{proof}

\smallskip

\begin{remark}
	A similar identity was proved in \cite{AGIS} by using integration over certain tropical projective spaces. The proof we provide here is purely elementary.
\end{remark}

\end{document}